\documentclass[runningheads,envcountsame]{llncs}
\usepackage[T1]{fontenc}
\usepackage{comment}

\usepackage{graphicx}
\usepackage{a4}
\usepackage[english]{babel}
\usepackage[dvipsnames,usenames]{color}
\usepackage{mathtools}
\usepackage{bbm}

\usepackage{verbatim,graphics,indentfirst} 
\usepackage{url}
\usepackage{stmaryrd}
\usepackage{booktabs}

\usepackage{latexsym}
\usepackage{enumerate}
\usepackage{graphicx}
\usepackage[T1]{fontenc}
\usepackage{lmodern}

\usepackage[utf8]{inputenc}
\usepackage{graphicx}
\usepackage{amsfonts}
\usepackage{nicefrac}
\usepackage{array}
\usepackage{appendix}

\usepackage{tikz}
\usepackage{tikz-qtree}
\usepackage[]{forest}

\usepackage{xspace}
\usepackage{marginnote}
\usepackage{proof}

\usepackage{pdflscape}

\usepackage{tikz}
\usepackage{tikz-cd}
\usetikzlibrary{tikzmark,decorations.pathreplacing,calc}

\usepackage{mathdots}
\usepackage[margin=1.5in]{geometry}

\makeatletter
\newcommand{\dotfrac}[2]{
\mathchoice
{\ooalign{$\genfrac{}{}{0pt}{0}{#1}{#2}$\cr\leavevmode\cleaders\hb@xt@ .22em{\hss $\displaystyle\cdot$\hss}\hfill\kern\z@\cr}}
{\ooalign{$\genfrac{}{}{0pt}{1}{#1}{#2}$\cr\leavevmode\cleaders\hb@xt@ .22em{\hss $\textstyle\cdot$\hss}\hfill\kern\z@\cr}}
{\ooalign{$\genfrac{}{}{0pt}{2}{#1}{#2}$\cr\leavevmode\cleaders\hb@xt@ .22em{\hss $\scriptstyle\cdot$\hss}\hfill\kern\z@\cr}}
{\ooalign{$\genfrac{}{}{0pt}{3}{#1}{#2}$\cr\leavevmode\cleaders\hb@xt@ .22em{\hss $\scriptscriptstyle\cdot$\hss}\hfill\kern\z@\cr}}
}
\makeatother

\renewcommand{\phi}{\varphi}
\newcommand{\A}{\Al}

\newcommand{\GenSubs}[1]{S}

\def\*#1{\mathbf{#1}}

\tikzset{
	treenode/.style = {align=center, inner sep=0pt, text centered},
	Ske/.style = {treenode, ellipse, double, draw=black,
		minimum width=6pt, thick},
	PIA/.style = {treenode, ellipse, black, draw=black,
		minimum width=6pt},
	Crit/.style = {treenode, rectangle, draw=black,
		minimum width=0.5em, minimum height=0.5em}
}

\newcommand{\tv}{\ensuremath{\mathbf 1}}

\newcommand{\OLor}{\lor_{\mathsf{OL}}}
\newcommand{\OLand}{\land_{\mathsf{OL}}}

\newcommand{\OLimp}{\ensuremath{\to_{\mathsf{OL}}}}
\newcommand{\OLneg}{\neg}

\newcommand{\la}{\langle}
\newcommand{\ra}{\rangle}
\newcommand{\alg}[1]{\mathbf{#1}}                  
\newcommand{\Al}[1][A]{\ensuremath{\mathbf{#1}} }
\newcommand{\nnot}{\mathop{\neg}}                  
\let\imp\Rightarrow                                
\newcommand{\pr}{\, \mathop{\preceq} \,}                 
 \newcommand{\eq}{\mathbin\equiv}

\newcommand{\uv}{\ensuremath{\top}}

\let\abs=\envert

\newcommand{\bv}{\ensuremath{\mathbf{0}}}

\newcolumntype{P}[1]{>{\centering\arraybackslash}p{#1}}

\newcommand{\uand}{\land_{\mathsf{K}}}
\newcommand{\uor}{\lor_{\mathsf{K}}}

\newcommand{\DFimp}{\to_{\mathsf{DF}}}

\newcommand{\Fimp}{\to_{\mathsf{F}}}

\renewenvironment{proof}{\noindent\emph{Proof.  }}{\hfill$\blacksquare$ \medskip}

\usepackage{amssymb}

\renewcommand{\la}{(}
\renewcommand{\ra}{)}

\newcommand{\sobocinski}{Soboci\'nski }

\begin{document}
%
\title{Indicative conditionals: from three to four values}
\author{Umberto Rivieccio\orcidID{0000-0003-1364-5003} \\
Miguel Muñoz Pérez\orcidID{0009-0001-3468-8208}}

%
\authorrunning{Rivieccio \& Muñoz}
%
\institute{Departamento de L\'ogica, Historia y Filosof\'ia de la Ciencia, \\
 UNED, Madrid, Spain \\
 \email{umberto@fsof.uned.es}\\
\email{mupemiguel99@gmail.com}
}

\maketitle              

\begin{abstract}
In this work we 
study 
how one may expand some three-valued logics of indicative conditionals by the addition of a fourth truth-value. This is achieved through the use of twist constructions, previously introduced in the literature in order to provide algebraic semantics for these logical systems. We first present some twist-representation results for the three-valued case, as well as for some relevant fragments and, then, we distinguish two ways of adding the new truth-value. We also give twist-representation results for these new cases. In the course of the paper we discuss the philosophical motivations behind the use of the new four-valued systems and some potential objections against these.

\end{abstract}

\section{Introduction}

The present  paper continues 
the study of logics of indicative conditionals 
that we initiated in the papers
\cite{wollic,tableaux}. The former of these  analyzes the algebraic semantics of several three-valued logical systems that have been proposed in the literature in order to provide sound accounts of implication as used in ordinary reasoning, and the latter deals with  axiomatizability issues of these systems. 
Our aim here is to describe some ways of extending logics of indicative conditionals, mainly through the use of \textit{twist-structures} (see, e.g. \cite{um20,FaOd23,Riv2x}).

\subsection{Background}

By \textit{indicative conditionals} we understand the simplest if-then statements employed in natural language, 
which are concerned with what could be true. They ought therefore to be distinguished from \textit{counterfactuals}, which deal with situations that are not possible from the outset. A classical example that illustrates such difference is provided by \emph{Adams' pair}: `if Lee Harvey-Oswald didn’t shoot JFK, someone else did' is an example of an indicative conditional, while `if Lee Harvey-Oswald had not shot JFK, someone else would have' is a counterfactual statement. One is clearly committed to the truth of the former sentence while, at the same time, may reject the latter (if, say, the speaker is not fond of theories involving the CIA or the KGB). In short, both kinds of conditional statements do not exclusively differ syntactically (say, in their grammatical mood): we may say that indicative conditionals have \emph{epistemological} flavor, while counterfactuals bear a more \emph{metaphysical} one~\cite{khoo}. 
\\

The literature on the possible formalizations of indicative conditionals is quite vast (see, e.g. \cite{egre20211,egre20212} for a survey). Here, we will exclusively deal with the trend that seeks to understand the indicative conditional as a truth-functional logical connective $\varphi \to \psi$, in contrast with the \textit{intensional} one, which gives up truth-functionality (see, e.g., \cite{sprenger}). In classical logic, $\varphi \to \psi$ is identified with the \emph{material conditional} and is thus deemed  equivalent to the (bivalent) disjunction $\neg \varphi \lor \psi$: 
`if $\varphi$ then $\psi$' is then false in case $\varphi$ is true and $\psi$ is false, and it is true otherwise. Leaving aside the famous paradoxes of material implication (see, e.g. \cite{egre20211}), several objections may be raised against this analysis, especially in cases in which the antecedent is false\footnote{Note how, in this very same way, the material conditional analysis trivializes counterfactual statements and, in particular, the distinction illustrated by Adam's pair.}. 

In these situations it seems that, while the material conditional is vacuously true, it fails to be informative. In order to verify the conditional statement $\varphi \to \psi$, one usually checks whether both $\varphi$ and $\psi$ hold while, when falsifying it, one expects to observe that $\varphi$ holds and that $\psi$ fails, so no actual information is obtained in case $\varphi$ fails. Hence, one may intuitively speak of a semantic \textit{gap}, i.e. one can regard such statements as lacking a definite truth-value. This approach can be traced back at least to Ramsey (1929), Reichenbach (1935), De Finetti (1936), Rhinelander and Quine (1950) (see, e.g., \cite{egre20211,egre20212,wollic}). Thus, within the truth-functional trend, we restrict ourselves to the \emph{suppositional} approach, condensed in the slogan: `uttering a conditional amounts to making a \emph{conditional assertion} [sic]: the speaker is committed to the truth of the consequent when the antecedent is true, but committed to neither truth nor falsity of the consequent when the antecedent is false' \cite{egre20211}.
\\

Among the possible ways of formalizing the above intuitions, if truth-functionality is to be preserved then one can simply add a new truth-value $\uv$ to the classical ones $\bv$ and $\tv$ in such a way that $\bv \to \psi = \uv$. It is also common to consider logics in which `non-falsity' is preserved under inferences (or, equivalently, where falsity is not introduced under inferences \cite{Ca08}), that is, logics having $\uv$ and $\tv$ as \emph{designated elements}\footnote{\emph{Prima facie}, it could be argued that one should naturally study truth-preserving logics instead. However, the formal behavior of these is somewhat unsatisfying: see \cite{wollic} for some notorious examples.}. Of course, the task of determining the behavior of $\uv$ with respect to the usual connectives is left open (intuitions may differ regarding, e.g., the value of $\uv \to \psi$). But these constraints are enough to determine the class of what we call (propositional) \textit{logics of indicative conditionals}.
\begin{figure}[h]
    \centering
    \begin{tabular}{@{}c|ccc@{}}
        \toprule
         $\OLand$ & \bv & \uv & \tv  \\
        \midrule
        \bv & \bv & \bv & \bv \\
         \uv & \bv  & \uv &  \tv\\
         \tv & \bv  & \tv &  \tv \\
        \bottomrule
    \end{tabular}
    \quad
    \begin{tabular}{@{}c|ccc@{}}
        \toprule
         $\OLor$ & \bv & \uv & \tv \\
        \midrule
        \bv & \bv & \bv & \tv \\
         \uv & \bv & \uv & \tv\\
         \tv & \tv & \tv & \tv \\
        \bottomrule
    \end{tabular}
    \quad
    \begin{tabular}{@{}c|ccc@{}}
        \toprule
         $\uand$ & \bv & \uv & \tv\\
        \midrule
         \bv & \bv & \bv & \bv\\
         \uv & \bv & \uv & \uv \\
         \tv & \bv & \uv & \tv\\
        \bottomrule
    \end{tabular}
    \quad
    \begin{tabular}{@{}c|ccc@{}}
        \toprule
         $\uor$ & \bv & \uv & \tv\\
        \midrule
         \bv & \bv & \uv & \tv\\
         \uv & \uv & \uv & \tv \\
         \tv & \tv & \tv & \tv\\
        \bottomrule
    \end{tabular}
    \\[1.2em]
    \begin{tabular}{@{}c|c@{}}
        \toprule
         & $\OLneg$ \\
         \midrule
         \bv & \tv\\
         \uv & \uv\\
         \tv & \bv\\
         \bottomrule
    \end{tabular}
    \quad
    \begin{tabular}{@{}c|ccc@{}}
        \toprule
         $\OLimp$ & \bv & \uv & \tv \\
        \midrule
        \bv & \uv & \uv & \uv\\
         \uv & \bv & \uv & \tv\\
         \tv & \bv & \uv & \tv\\
        \bottomrule
    \end{tabular}
    \quad
    \begin{tabular}{@{}c|ccc@{}}
        \toprule
         $\DFimp$ & \bv & \uv & \tv\\
        \midrule
         \bv & \uv & \uv & \uv\\
         \uv & \uv & \uv & \uv \\
         \tv & \bv & \uv & \tv\\
        \bottomrule
    \end{tabular}
    \quad
    \begin{tabular}{@{}c|ccc@{}}
        \toprule
         $\Fimp$ & \bv & \uv & \tv\\
        \midrule
         \bv & \uv & \uv & \uv\\
         \uv & \bv & \uv & \uv \\
         \tv & \bv & \uv & \tv\\
        \bottomrule
    \end{tabular}
\caption{Tables of three-valued connectives for indicative conditionals.}
\label{fig:truth-tables2}
\end{figure}

The distinctive formal behavior of logics of indicative conditionals can be motivated by appealing to some pragmatic considerations such as, v.g., verifying certain \emph{connexive principles} that are, in general, classical contingencies\footnote{Two notorious examples of these principles are given by \textit{Aristotle's} and \textit{Boethius' theses}, respectively: (A1) $\neg(\varphi \to \neg \varphi)$, (A2) $\neg(\neg \varphi \to \varphi)$ and (B1) $(\varphi \to \psi) \to \neg (\varphi \to \neg \psi)$, (B2) $(\varphi \to \neg \psi) \to \neg (\varphi \to \psi)$ (see, e.g., \cite{FaOd23}). For more details, see \cite{wollic}.} or other statements that make them \emph{contra-classical} and hence incomparable with classical logic (rather than sub-classical) \cite{wollic}. Alternatively, one may try to provide some pre-formal insight consistent with the philosophical considerations above, namely, considering $\to$ as behaving like a restriction of the material conditional that simply excludes the preceding problematic cases (see \cite[\S7.19]{humberstone}): one selects some truth (resp. falsity) conditions of classical implication and then defines them to be the truth (resp. falsity) conditions of $\to$ while leaving the remaining ones to be $\uv$. In this way, for example, $\varphi \to_\mathsf{DF} \psi$ is true \emph{just in case} both $\varphi$ and $\psi$ are true, and false \emph{just in case} $\varphi$ is true and $\psi$ is false; otherwise, this formula receives the value $\uv$. The truth conditions for three-valued conjunctions can be similarly defined in terms of those for classical conjunction (see below).
\\

The basic three-valued systems 
we will consider here are the following ones\footnote{
Note that our strategy above allows us to motivate \emph{more} logics than the presented here. Part of our future work will consist in studying these logics.}
. Recall that a common way of defining a logical system is to consider the \emph{logic induced by a matrix} $(\*{A}, F)$ which is formed by an algebra $\*{A}$ (roughly, comprising the corresponding truth tables) and a subset $F \subseteq A$ of designated values (i.e., those to be preserved under inferences). Let us fix $A_3 := \{\bv, \uv, \tv\}$. 
\begin{itemize}
\setlength{\itemsep}{9pt}
    \item[1.] \textit{De Finetti's logic} DF is 
    induced by $(\mathbf{DF_3}, \{\uv, \tv\})$, where $\mathbf{DF_3} := (A_3; \neg, \land_\mathsf{K}, \lor_\mathsf{K}, \to_\mathsf{DF})$ \cite{deFin,egre20211,egre20212}. 
    \item[2.] \textit{Cooper's logic} OL is 
    induced by $(\mathbf{OL_3}, \{\uv, \tv\})$, where $\mathbf{OL_3} := (A_3; \neg, \land_\mathsf{OL}, \lor_\mathsf{OL}, \to_\mathsf{OL})$ \cite{Cooper1968,Riv2x}. 
    \item[3.] \textit{Cantwell's logic of conditional negation} CN is 
    induced by $(\mathbf{CN_3}, \{\uv, \tv\})$, where $\mathbf{CN_3} \newline := (A_3; \neg, \land_\mathsf{K}, \lor_\mathsf{K},\to_\mathsf{OL})$ \cite{Ca08}. This system, as shown in~\cite{wollic}, is definitionally equivalent to \textit{Farrell's logic} F, that is, the consequence induced by $(\mathbf{F_3}, \{\uv, \tv\})$, with $\mathbf{F_3} := (A_3; \neg, \land_\mathsf{K}, \lor_\mathsf{K}, \to_\mathsf{F})$ \cite{Fa79}. 
\end{itemize}


\subsection{Algebraic preliminaries}

In this subsection we 
introduce the key notions that we will use in the what follows. 
For more details on the basic notions of algebraic logic see, e.g., \cite{Um17,FJa09,F16}; for undefined universal algebra notions see, e.g., \cite{BuSa00}. 
The main idea behind the use of twist-structures consists in representing the algebraic semantics associated with a given logic (in our case, DF, CN and OL) as a special kind of power construction in terms of more standard algebraic structures (see, e.g. \cite{Riv14,um20}). We begin by introducing these well-known classes of algebras:
\begin{definition}
\label{demorgan}
    A distributive lattice $(L; \land, \lor)$ is a \emph{De Morgan lattice} if there is an operation $\neg: L \to L$ such that both $\neg\neg a = a$ and $\neg(a \land b) = \neg a \lor \neg b$ hold for every $a, b\in L$. A De Morgan lattice $L$ is called \emph{centered} if there is some element $c \in L$ such that $\neg c = c$. If $L$ has two distinct centers, we say that it is \emph{bi-centered}. In case a De Morgan lattice is bounded, we speak of a \emph{De Morgan algebra}. 
\end{definition}

We will interested in the following kind of De Morgan algebras \cite{Cig86}: 

\begin{definition}
\label{kleene}
    A \emph{Kleene algebra} is a De Morgan lattice $L$ in which it holds that $\neg a \land a \leqslant \neg b \lor b$, for every $a, b \in L$. Note that, in case a Kleene algebra has a center, such center is unique (otherwise, given two different centers $c, d$ we would have that $c = c \land \neg c \leqslant d \land \neg d = d$ and, by symmetry, that $c = d$).
\end{definition}


Let us also recall that: 
\begin{definition}
\label{boolean}
   A \emph{generalized Boolean algebra} $(B; \land, \to)$ is the $0$-free reduct of a Boolean algebra $(B; \land, \to, 0)$. Equivalently, one can define it as an algebra $(B; \land, \lor, \to, 1)$ verifying that $(B; \land, \lor, 1)$ is an upward bounded lattice, that $\land, \to$ form a \emph{residuated pair} (i.e. $x \land y \leqslant z$ iff $y \leqslant x \to z$) and that \emph{Peirce's law} (i.e. $(x \to y) \to x = x$) holds. 
\end{definition}

It will be useful to remember that: 

\begin{definition}
\label{complemented}
    Given a bounded lattice $\*L := (L; \land, \lor, 0, 1)$, we say that $a \in L$ is \emph{complemented} in case there is some $b \in L$ such that $a \land b = 0$ and $a \lor b = 1$. We say that $\*L$ is \emph{complemented} in case each element of $L$ is complemented. Note that a complemented De Morgan algebra is a Boolean algebra. 
\end{definition}

We will also make use of the following definition (see \cite{celani08} for more details): 
\begin{definition}
\label{tarski}
A \emph{Tarski algebra} is an algebra $\mathbf{T}  := (T ;\to, 1)$ of type $(2, 0)$ that satisfies the following quasi-equations: (i) $x \to( y \to x) = 1 $, (ii) $x \to (y \to z) = (x \to y) \to (x \to z)$, (iii) $x \to y = y \to x = 1$ implies $x = y$ and (iv) $((x \to y) \to x ) \to x = 1$. In fact, Tarski algebras form a variety and characterize precisely the $\{ \to, 1 \}$-fragment (equivalent to the $\{ \to, \lor, 1 \}$-fragment) of Boolean algebras. 
\end{definition}

\paragraph{Notation.} In general, we will distinguish the logical connectives (defined by a concrete truth table) from the algebraic operations. The former will usually be accompanied by a subindex or some other typographical indication, while the latter (as in the previous definitions) will be written without any  index. 
{For instance, we will denote the three-valued conjunction and the  quasi-conjunction~\cite{egre20211}, respectively, by $\land_\mathsf{K}$ and $\land_\mathsf{OL}$, whereas their algebraic counterparts will appear as $\land$ and $*$, respectively. The same goes with $\lor$, $+$ etc.}

\section{The three-valued case}

In this section we illustrate how to represent the 
algebraic models of our logics in terms of twist-algebras. The overall strategy is common through all the examples considered. First, one identifies a class of algebras
that provides a semantics
for the logic at hand; 
then, one defines a class of twist-structures 
that belong to this class,
and, finally, one proves a \textit{twist representation} result, that is, that each member of the abstractly defined class is isomorphic to a twist-algebra.


\subsection{Cooper} 

We begin by reviewing the case of Cooper's Logic of Ordinary Discourse (OL), which was the first among the logics of indicative conditionals to be studied from this standpoint (see \cite{Riv2x,Greati23}). 
Cooper's original motivations for OL \cite{Cooper1968} were quite rich, his main point being that the logic of natural language radically diverges from classical logic and, in particular, the material conditional. He intended to prove this claim by exhibiting an exhaustive list of inferences sanctioned by informal speech but falsified by classical logic and, conversely, verified in the latter while falsified in the former. Following our comments from the introduction, one can motivate Cooper's connectives as follows. Regarding the conditional:  
$$\varphi \to_\mathsf{OL} \psi = \tv \text{ iff } \varphi \neq \bv \text{ and } \psi = \tv, \quad \varphi \to_\mathsf{OL} \psi = \bv \text{ iff } \varphi \neq \bv \text{ and } \psi = \bv.$$
Similarly, for conjunction one has: 
$$\varphi \land_\mathsf{OL} \psi = \tv \text{ iff } (\varphi = \tv \text{ and } \psi \neq \bv) \text{ or } (\varphi \neq \bv \text{ and } \psi = \tv), \quad \varphi \land_\mathsf{OL} \psi = \bv \text{ iff } \varphi = \bv \text{ or } \psi = \bv.$$
It is also worth noting that Cooper's original proposal did not allow atomic propositions to receive truth-value $\uv$, which is consistent with the idea that $\uv$ arises in conditional contexts. However, the resulting system would fail to be structural\footnote{See \cite{Riv2x} for more technical details on this. In~\cite{egreprob}, however, \emph{Cooper valuations} are employed and motivated as having some theoretical appeal.}. This is why we restrict ourselves to the structural version of OL (denoted $\mathrm{sOL}$ in \cite{Riv2x}), namely, the logic induced the matrix $(\mathbf{OL_3}, \{\uv, \tv\})$ (see above).

\paragraph{Notation.} In this section (and, unless stated otherwise, throughout the whole paper) we shall use the following abbreviations. Let $\abs{x} := x \to x$, and abbreviate $x \to y = \abs{x \to y}$ as $x \pr y$.
Further, we write $x \eq y$ when $x \pr y \pr x$. We view these both as terms/equations
and as relations defined on a given algebra. In the cases of interest,
the relation $ \pr $ will in general be a pre-order (i.e.~reflexive and transitive)
and $\eq$ will be an equivalence relation compatible with some (but not all) the algebraic operations.
In the case of OL-algebras (see below), the relation $\eq$ can also be defined via the term
$\Diamond x : = \nnot x \to x$, for one has $x \eq y$ if and only if $\Diamond x = \Diamond y$.
\\

It is known that OL is algebraizable~\cite[Thm.~3.1]{Riv2x}, and 
the corresponding equivalent semantics is
the variety of \emph{OL-algebras}~\cite[Def. 2.1]{Riv2x}: 
\begin{definition}
\label{def:oalg}
An \emph{OL-algebra} is an algebra $\Al = (A ; *, \to, \nnot)$ of type
$(2, 2, 1)$
such that
$(A ;  *  )$ is a semilattice (with order $\leq$) and
the following equations are satisfied:

\begin{enumerate}[(OL1)
]
\item \label{Itm:OL1eq} $\abs{x}  \to y = y $,
 \item \label{Itm:OL2eq} $x \to (y \to z) = (x \to y) \to (x \to z)$,
\item \label{Itm:OL3eq} $((x \to y ) \to x) \to x = \abs{x}$,
\item \label{Itm:OL4eq} $x  *  y = x  *  (x \to y) $,
 \item \label{Itm:OL5eq} $x \to (y \to z) = (x  *  y) \to z$,
\item \label{Itm:OL6eq} $x \to (y  *  z) = (x \to y)  *  (x \to z)$,
\item  \label{Itm:OL7eq} $\neg \neg x = x$,
\item \label{Itm:OL8eq} 
$ x \to y = (\nnot x \to x ) \to y $,
\item \label{Itm:OL9eq} $ x \to \nnot y = \nnot (x \to  y) $,
\item \label{Itm:OL11eq} $ \nnot (x  *  y)  \pr  (x \to \nnot y)  *  (y \to \nnot x) $,
\item \label{Itm:OL13eq} $ (x \to \nnot y)  *  (y \to \nnot x)   \pr  \nnot (x  *  y)  $,
\item \label{Itm:OL10eq} 
$x = \abs{x}  *  (\nnot x \to x) $,
\item \label{Itm:OL12eq} 
$\nnot x \to x \leq (x \to  y ) \to (\nnot y \to y )$.
\end{enumerate}
\end{definition}
On every OL-algebra, a second semilattice operation $\lor$ may be obtained as the De Morgan dual of $ * $ (as in Definition~\ref{OLtwist} below).
Note, however, that in general the absorption law between $ * $ and $\lor$ fails to hold,
so the two operations determine two independent semilattice structures, not one lattice. See, e.g.,~\cite{Riv2x} for further details on other term-definable (semi-)lattice orders.
We present below the twist construction for OL-algebras~\cite[Def.~2.3]{Riv2x}:

\begin{definition}
\label{OLtwist}
     Let $\*B := (B; \land, \to)$ be a generalized Boolean algebra (i.e.~the $0$-free subreduct of a Boolean algebra). The \emph{OL-twist-algebra over $\*B$} is the algebra $\mathbf{B}^{\bowtie}$ having as universe $\{(a_1, a_2)\in B^2 \, | \, a_1 \to a_2 = a_2\}$\footnote{Equivalently, one may require that $a_1 \lor a_2 = 1$ instead (recall that $x \lor y = (x \to y) \to y$ when $\to$ is the Boolean implication).} and, as operations:
\begin{itemize}
    \setlength \itemsep{0pt}   
    \item[i.] $(x_1, x_2) * (y_1, y_2) := (x_1 \land y_1, (x_1 \to y_2)\land (y_1 \to x_2)),$ 
    \item[ii.] $(x_1, x_2) \to (y_1, y_2) := (x_1 \to y_1, x_1 \to y_2),$ 
    \item[iii.] $\neg (x_1, x_2) := (x_2, x_1)$.
\end{itemize}
The operation $+$ is obtained as the De Morgan dual of $*$, namely,
$x + y : = \nnot (\nnot x * \nnot y)$.
\end{definition}

The usual twist representation result holds in this case: every OL-twist-algebra constructed in the previous way is in fact an OL-algebra \cite[Prop. 2.4]{Riv2x} and, 
moreover, every OL-algebra is of this form
\cite[Thm. 2.5]{Riv2x}. This fact is quite useful, and allows us 
to prove, e.g., that the class of OL-algebras 
is generated as a variety by
$\mathbf{OL_3}$ \cite[Cor. 2.7]{Riv2x}. 
The strategy followed in \cite{Riv2x} in order arrive at the desired twist representation is that of defining $\Diamond x := \neg x \to x$ and $x \eq y$ by $\Diamond x = \Diamond y$ \cite[Prop. 2.2(i)]{Riv2x}, then checking that $\eq$ is a congruence on the $\{\land, \to \}$-reduct and that the corresponding quotient algebra is a generalized Boolean algebra~\cite[Prop. 2.2(vi)]{Riv2x}. 
Equivalently, instead of using the quotient $\*A / \eq$, one can obtain the factor algebra also
as the direct image of the $\Diamond$ operator, that is,
the algebra $\Diamond(\*A) := (\Diamond[A]; \land, \to)$ whose operations
are just the restrictions of those on 
$\*A$. One thus arrives at the following representation result \cite[Thm. 2.5]{Riv2x}:
\begin{theorem}[OL-twist representation]
\label{OLtwistrep}   
    Every OL-algebra $\*A$ is isomorphic to a OL-twist-algebra over $\Diamond(\*A)$ via the map  $\iota \colon \*A \to \Diamond (\*A) \times \Diamond (\*A): a \mapsto (\Diamond a, \Diamond \nnot a)$.
\end{theorem}
Note how further properties of OL-algebras may be needed through the proof in order for $\iota$ to be an embedding. We redirect the reader interested in the details to \cite{Riv2x}.

\subsection{Fragments of OL}
The logical language of OL $\{ \land, \to, \nnot \} $
is quite expressive, for the three-element algebra
that generates the variety of OL-algebras is quasi-primal~\cite[Thm.~12]{Riv2x}, and it becomes primal
(so all three-valued functions become expressible)
if we further add the constant $\tv$ (or, equivalently, $\bv$).
On the other hand, by simply deleting the negation from the  set $\{ \land, \to, \nnot \} $ one obtains fragments that coincide with
the classical ones. This suggests as objects of obvious interest the two fragments 
$\{\to, \nnot \}$ and $\{\nnot, \land \}$, which are both algebraizable and non-classical;
 other term-definable connectives in OL (e.g.~the weak Nelson implication)
  may also be employed 
in establishing algebraizability\footnote{The `term-definable fragment' of OL corresponding to the primitive negation and the definable weak Nelson implication is studied in detail in~\cite{Riv14}. 
 }. The fragment $\{\nnot, \land \}$, which is equivalent to $\{\nnot, \imp \}$, corresponds precisely to
\emph{Soboci\'nski's logic}~\cite{sobocinski,Pa72}, while the others, as far as we are aware, have not yet been studied 
in isolation. In this section we shall take a look at these two algebraizable sub-logics of OL,
leaving for future work the study of other fragments that may be of independent interest,
such as $\{ \land, \lor \}$ or $\{ \imp \}$.

\subsection*{The $\{\to, \nnot\}$-fragment}

One of the advantages of twist representations is that they can be easily adapted to characterize fragments of
logics (resp.~algebras) as long as the connectives on which the representation hinges are retained
(see e.g.~\cite{JaRi21}). Here we  
illustrate this 
by focusing the $\{\to, \nnot\}$-fragment of OL, which we believe to be of particular logical interest: indeed, it is precisely this fragment  -- which, we note, OL shares with Farrell's logic -- 
that makes OL a connexive and a contra-classical logic (by contrast, the $\{\land, \to\}$-fragment is classical,
and $\{\land, \lor, \nnot \}$-fragment is sub-classical). We shall focus mainly on algebraic models, leaving the issue of axiomatizing the corresponding logic to a future publication (but see~\cite{Greati23}); 
this can be easily done by using the algebraizability transformer to translate back the equations that 
form a basis for the corresponding variety (Definition~\ref{d:inOL} below).
\\


The twist construction and representation follow closely those introduced above for the algebras
in the full language, with two important differences:
(i) as one would expect, the factor algebras involved correspond to a fragment (the purely implicational) of Boolean algebras;
(ii) because of the language restriction, in this case we will only be able to obtain an embedding (not an isomorphism)
result.

\begin{definition}
\label{d:OLnegimptwist}
Let $\mathbf{T} := (T ; 
\to ,  
1
)$
be a Tarski algebra.
An \emph{inOL-twist-algebra} is an algebra
$ 
\mathbf{T}^{\bowtie}$
with universe 
$\{ \la a_1 , a_2  \ra \in T \times T : a_1 \to a_2 = a_2 
\}
$
and 
operations given,
by:
\begin{enumerate}[i.]
    \item $\nnot ( x_1 , x_2 ) := (x_2,   x_1)$,
    \item $(x_1  , x_2) \to ( y_1 , y_2 )  : = (x_1  \to  y_1  ,   x_1 \to  y_2).$
\end{enumerate}
%
An \emph{
inOL-twist-algebra over $\mathbf{T}$}
is any subalgebra $\mathbf{A} \leqslant \mathbf{T}^{\bowtie} 
$
satisfying
$\pi_1[A] = T$. 
\end{definition}

On every Tarski algebra, 
joins exist and are given by $x \lor y := (x \to y) \to y$ (see, e.g., \cite{celani08}), so one may equivalently rewrite the condition
$a_1 \to a_2 = a_2 $ as $a_1 \lor a_2 = 1$.
It is easy to verify that the set $\mathbf{T}^{\bowtie} $ is itself closed under the operations
(indeed, this follows from the corresponding result on OL-algebras).
\\

We abbreviate $\abs{x} := x \to x$ and we write $x \pr y$ for $x \to y = \abs{x \to y}$
and $x \eq y$ when $x \pr y \pr x$.
We shall also employ the following term (cf.~\cite[p. 279]{JaRi21}): 
$$
q(x,y,z) := (x \to y) \to ((y \to x) \to ((\nnot x \to \nnot y) \to ((\nnot y \to \nnot x) \to z))).
$$



We next introduce our candidate for an abstract definition of the class of $\{\to, \nnot\}$-subreducts of OL-algebras:
note that this is, by definition, a variety. 

\begin{definition}
\label{d:inOL}
An \emph{implication-negation OL-algebra} (for short, \emph{inOL-algebra}) is an algebra $\mathbf{A} := (A ; 
 \to, \nnot 
)$ 
of type $(2, 
1 
)$ satisfying the following properties:
 
\begin{enumerate}[({in}OL1)
]
\item \label{Itm:inOL1eq} $\abs{x}  \to y = y $,
 \item \label{Itm:inOL2eq} $x \to (y \to z) = (x \to y) \to (x \to z) =  y \to (x \to z)$,
\item \label{Itm:inOL3eq} $((x \to y ) \to x) \to x = \abs{x}$,
\item \label{Itm:inOL4eq} $x \pr y \to x$,
\item  \label{Itm:inOL7eq} $\nnot \nnot x = x$,
\item \label{Itm:inOL8eq} 
$ x \to y = (\nnot x \to x ) \to y $,
\item \label{Itm:inOL9eq} $ x \to \nnot y = \nnot (x \to  y) $,

\item \label{Itm:inOL10eq} $q(x,y,x) = q(x,y,y) $.
\end{enumerate}
\end{definition}

The above-defined class is therefore a variety; note that the only equations that do not appear in
the definition of OL-algebras (Definition~\ref{def:oalg}) are the second one of~(inOL\ref{Itm:inOL2eq})
and the last,~(inOL\ref{Itm:inOL10eq}). 
 Now, we can consider $\pr$ as a relation on each inOL-algebra $A$, as well as its 
symmetrization $\eq$ defined by $x \eq y$ whenever $x \pr y \pr x$ (see the Notation paragraph at the beginning of this section). 

\begin{lemma}
\label{lem:inolfirsancheck}
Let $(A; \to, \neg)$ be a inOL-algebra. Then, the following holds: 
\begin{enumerate}[(i)]

\item $\pr$ is a preorder relation (i.e. reflexive and transitive),

\item $x \pr y$ implies 
$y \to z \pr x \to z$ and $z \to x \pr z \to y$,

 \item  In consequence, $\eq$ is a congruence  on the reduct $( A; \to )$,

\item 
$(A; \to) /  \eq $ is 
a Tarski algebra with top element $(a \to a) / \eq$ (for an arbitrary $a \in A$), 

\item $x \eq y$ and $\nnot x \eq \nnot y$ 
imply $x =y$,

\end{enumerate}
\end{lemma}

\begin{proof}
(i). 
For each $a \in A$, we have $a \pr a$
by~(inOL\ref{Itm:inOL1eq}). 
To verify that $\pr$ is transitive, assume $a \pr b$ and $b \pr c$.
Using~(inOL\ref{Itm:inOL2eq}), 
we  have
$a \to (b \to c) = (a \to b) \to (a \to c) = \abs{a \to b} \to  (a \to c) = a \to c$.
On the other hand, from $b \pr c$ 
we have
$a \to (b \to c) = a \to \abs{b \to c} = (b \to c) \to (a \to (b \to c)) = \abs{(b \to c) \to (a \to (b \to c))}$.
By~(inOL\ref{Itm:inOL1eq}),
 we conclude that $(a \to c) \to d = d$ for all $d \in A$. In particular
$a \to c = \abs{a \to c}$, whence $a \pr $c, as required.

(ii).
Assume $a \pr b$. 
%
Then, using~(inOL\ref{Itm:inOL2eq}),
we have $c \to (a \to b) = c \to \abs{a \to b }  = (a \to b) \to (c \to (a \to b)) =
\abs{c \to (a \to b)}$. But since $ c \to (a \to b) = (c \to a) \to (c \to b)$, we conclude that $c \to a \pr c \to b$, as required.
Further observe that the premise $a \pr b$, 
by~(inOL\ref{Itm:inOL1eq}),
gives us
$(a \to b) \to ((b \to c ) \to (a \to c)) = (b \to c ) \to (a \to c)  $.
On the other hand,  using~(inOL\ref{Itm:inOL2eq}),
we have
$(a \to b) \to ((b \to c ) \to (a \to c)) 
= (b \to c ) \to ((a \to b)  \to (a \to c)) 
= (b \to c ) \to (a \to (b  \to c)) 
= a  \to  \abs{b  \to c}$.
Thus, 
$(b \to c ) \to (a \to c)  = a  \to  \abs{b  \to c}$, which gives us
$b \to c \pr a \to c$, as desired.

(iii). Item (i) above implies that 
$\eq$ is an equivalence relation, which is, by item (ii), 
compatible with 
$\to$.

(iv). First observe $[a \to a ]$ is the top element of the Tarski algebra, that is,
$[b] \to [a \to a ] = [a \to a ]$, for all $a,b \in A$. Indeed, we have
$b \to (a \to a ) = (b \to a) \to (b \to a)$, entailing that both
$b \to (a \to a ) \pr a \to a$ and $a \to a \pr b \to (a \to a )$.
Hence, we may denote $1 := [a \to a]$.
Further, let us note that, for all $a,b \in A$, we have
$a \pr b$
iff 
$a \to b = (a \to b) \to (a \to b)$
iff
$[a \to b] = [(a \to b) \to (a \to b)]$ 
iff $[a] \leq [b]$.
Using these observations, it is easy to realize that equations (ii) and (iv) in the 
definition of Tarski algebras are enforced, respectively, by~(inOL\ref{Itm:inOL2eq}) and~(inOL\ref{Itm:inOL3eq}),
while (iii) is obviously satisfied. As to (i), it suffices to compute
$a \to (b \to a ) = b \to (a \to a ) = (b \to a) \to (b \to a)$.

(v). Assuming $a \eq b$ and $\nnot a \eq \nnot b$ and using~(inOL\ref{Itm:inOL1eq}),
it is easy to verify that $ q(a,b,a)= a$ and $q(a,b,b) = b$. The required result then follows from~(inOL\ref{Itm:inOL10eq}).
\end{proof}

\begin{theorem}[inOL-twist representation]
\label{thm:inoltwrep}
Every inOL-algebra $\Al$ embeds into the inOL-twist-algebra
$[(A; \to) / \eq]^{\bowtie}$.
\end{theorem}
\begin{proof}
Lemma \ref{lem:inolfirsancheck} ensures that we have the correct ingredients to build a twist-algebra. 
The embedding $\iota \colon \mathbf{A} \to \mathbf{A}/\eq_1 \times \mathbf{A}/\eq_2$ is defined, as expected, by
$\iota(a) := ([a], [\nnot a])$, which
is injective by  item (v) of Lemma \ref{lem:inolfirsancheck}.
Let us check that $\iota$ preserves the algebraic operations. 
The case of the negation 
is immediate. As for the implication, we have:
\begin{align*}
\iota (a \to b) 
& = ([a \to b], [ \nnot (a \to b)]) \\
& =  ([a \to b ], [a \to \nnot b]) & \text{by~(inOL\ref{Itm:inOL9eq})}\\
& =  ([a] \to [b ], [a] \to [\nnot b]) \\
& =  ([a], [\nnot a ]) \to  ([b], [\nnot b]) \\
& = \iota (a) \to \iota (b).
\end{align*}
Finally, observe that,
for every element $([a], [\nnot a]) \in \iota[A]$, we have
$[a ] \to [\nnot a ] = [ \nnot a ] $, as required by Definition~\ref{d:OLnegimptwist}.
Indeed, $\nnot a \pr a \to \nnot a $ holds by~(inOL\ref{Itm:inOL4eq});
and  $a \to \nnot a \pr \nnot a $ by~(inOL\ref{Itm:inOL8eq}), for we have
$(a \to \nnot a) \to \nnot a = (\nnot \nnot a \to \nnot a) \to \nnot a =  \nnot a \to \nnot a$.
\end{proof}

Theorem~\ref{thm:inoltwrep} allows us to identify any inOL-algebra $\mathbf{A}$ with a inOL-twist-algebra
$\mathbf{A} \leqslant \mathbf{T}^{\bowtie}$. This greatly simplifies calculations, and in particular it makes it easy to 
recover a number of properties
that are known to hold on OL-algebras (see \cite[Prop. 2.2, \S4]{Riv2x}).
Let us state a few.

\begin{proposition}
\label{prop:inOLprops}
 Given an inOL-algebra $\mathbf{A}$, the following holds:
\begin{enumerate}[(i)]

\item $x \to (\nnot x \to x ) $ defines a constant (that we denote, as before, by $\uv$),
\item $x \pr y $ \ iff \ $x \to \Diamond y = \uv$ \ iff \ $\Diamond x \to \Diamond y = \uv$,

\item $x \eq y $ \ iff \ $\Diamond x = \Diamond y $,
\item $x = y $ \ iff \ $\Diamond x = \Diamond y $ and $\Diamond \nnot x = \Diamond \nnot y $,

\item The $\land_\mathsf{K}$ partial order $\leqslant_\mathsf{K}$ on $\mathbf{A}$ is defined by 
\ $x \leqslant_\mathsf{K} y$ \ iff \ 
[$x \pr y$ and $\nnot y \pr \nnot x$],

\item The $\land_\mathsf{OL}$ partial order $\leqslant_\mathsf{OL}$ on $\mathbf{A}$ is defined by 
\ $x \leqslant_\mathsf{OL} y$ \ iff \ 
[$x \pr y$ and $x \pr \nnot x \to \nnot y \pr (\nnot x \to y) \to y$].
\end{enumerate}
    \end{proposition}
Regarding the defined operator $\Diamond$ that appears on the preceding Proposition,
it is instructive to observe that, on twist-algebras, one has 
$\Diamond \la a_1, a_2 \ra = \la a_2 \to a_1, a_1 \to a_1 \ra = \la a_1, 1 \ra $.
Thus $\Diamond$ has the important effect of keeping the first component of each pair while 
obliterating the second one: hence $\Diamond$ preserves the implication operator, 
and 
the forward image $\Diamond [A]$ is the universe of a Tarski
subalgebra of the twist-algebra $\A$.
We will see in the next sections that a similar operator can be defined
via different terms, according to  the algebraic language that is available (see e.g.~Subsection~\ref{ss:df}).

A twist representation often allows us to give a description of congruences
on $\mathbf{A} \leqslant \mathbf{T}^{\bowtie}$ in terms of those of 
 $\Al[T]$. Such a result was established for OL-algebras in~\cite[Thm.~4]{Riv2x}; we proceed 
 to prove its analogue for inOL-algebras.

\begin{lemma}
    \label{l:inOLcong}
    Let $\A$ be an inOL-algebra and $\theta$ a congruence on $\A$. Then, for all $a, b \in A$,
    $$
    \la a, b \ra \in \theta \ \  \text{ if, and only if, } \ \  \la \Diamond a, \Diamond b \ra, \, \la \Diamond \nnot a, \Diamond \nnot b \ra\in \theta 
    $$
\end{lemma}
\begin{proof}
   The ``only if'' part is straightforward. For the converse,
    assume $\la \Diamond a, \Diamond b \ra, \la \Diamond \nnot a, \Diamond \nnot b \ra\in \theta $.
    Then, in the quotient $\A  / \theta$ (which also belongs to the variety), we have $\Diamond [a] = [\Diamond a] = [\Diamond b] = \Diamond [a]$
    and, similarly, $\Diamond \nnot [a] = \Diamond \nnot [b]$. Then, by Proposition~\ref{prop:inOLprops} (iv),
    we conclude $[a] = [b]$, as required. 
\end{proof}

\begin{proposition}
    \label{p:inOLsi}
    If $\mathbf{A} \leqslant \mathbf{T}^{\bowtie}$ is a subdirectly irreducible inOL-algebra,
    then $\Al[T]$ is a subdirectly irreducible Tarski algebra.
\end{proposition}
\begin{proof}
To simplify the computations, let $T = \Diamond [A]$. Let $\theta$ be the monolith congruence on
$\A$. We claim that $\theta' : = \theta \cap (T \times T)$ is the monolith on $\Al[T]$.
It is clear that $\theta'$ is a congruence on $\Al[T]$.
To show that $\theta'$ is non-trivial, let $\la a,b \ra \in \theta$ be such that $a \neq b$. Then
$\la \Diamond a, \Diamond b \ra, \, \la \Diamond \nnot a, \Diamond \nnot b \ra\in \theta'$ and, by Proposition~\ref{prop:inOLprops} (iv),
either $ \Diamond a \neq \Diamond b$ or $\Diamond \nnot a \neq \Diamond \nnot b$.
In either case we conclude that 
$\theta' \neq Id_{T}$. To conclude the proof,
assuming $\eta'$ is a non-trivial congruence on $\Al[T]$, we verify that $\theta' \subseteq \eta'$.
To establish this, it suffices to show that $\eta' = \eta \cap (T \times T)$ for some non-trivial
congruence $\eta$ on $\A$. Define $\eta : = \{ \la a,b \ra \in A \times A : \la \Diamond a, \Diamond b \ra, \, \la \Diamond \nnot a, \Diamond \nnot b \ra\in \eta'\}$. It is easy to verify that $\eta$ is actually a congruence on
$\A$ (to prove compatibility with $\to$
one can use the fact that, for all $c \in A$, we have
$\Diamond \nnot (a \to c) = \Diamond (a \to \nnot c) = \Diamond a \to \Diamond \nnot c$ and, similarly,
$\Diamond \nnot (c \to a) = \Diamond (c \to \nnot a) = \Diamond c \to \Diamond \nnot a$).
Since $\eta'$ is non-trivial, let $\la a,b \ra \in \eta'$ be such that $a \neq b$.
Since $\Diamond x = \Diamond \Diamond x$ holds generally, the assumption that $a,b \in T$ gives us
$\Diamond a = a$ and $\Diamond b = b$. Note that 
$\Diamond \nnot \Diamond x = \uv$ is also a valid equation on inOL-algebras.
Thus, we have
$\la \Diamond a, \Diamond b \ra = \la a, b \ra \in \eta'$
and $ \la \Diamond \nnot a, \Diamond \nnot b \ra 
= \la \Diamond \nnot \Diamond a, \Diamond \nnot \Diamond b \ra
= \la \uv, \uv \ra \in \eta'$, which together give us 
$\la a,b \ra \in \eta$. So $\eta$ is non-trivial, whence we conclude that
$\theta \subseteq \eta$ and $\theta' \subseteq \eta'$. 
\end{proof}

It is easy to see that the proof of Proposition~\ref{p:inOLsi}
actually shows that the congruence lattice of each $\mathbf{A} \leqslant \mathbf{T}^{\bowtie}$
is isomorphic to that of the underlying Tarski algebra $\mathbf{T}^{\bowtie}$, thereby
implying, for instance, that inOL-algebras form a congruence-distributive variety.
More directly, the lemma tells us that
the only subdirectly irreducible inOL-algebra is the 
$\{\to, \nnot\}$-reduct of the three-element OL-algebra,
which embeds into the only subdirectly irreducible two-element Tarski algebra.
Hence, we have the following (for finitely-generated quasi-varieties, see e.g.~\cite[Thm. 3.6(ii)]{ClDa98}).

\begin{corollary}
    \label{c:inOLvareen}
    The class of inOL-algebras is generated, both as a variety and as a quasi-variety, 
    by its three-element member. Hence,
    inOL-algebras is precisely the class of $\{\to, \nnot\}$-subreducts of OL-algebras. 
\end{corollary}

\subsection*{The $\{ \imp, \nnot\}$-fragment, or Soboci\'nski's logic}

In this subsection we sketch a twist-algebra characterization for the algebraic models of
the $\{ \land, \lor, \nnot\}$-fragment of OL, which is equivalent to the $\{ \imp, \nnot\}$-fragment.
The truth table of $\imp$ on the three-element algebra is displayed in Figure~\ref{fig:truth-tables-deriv},
which also shows the connective $\supset$ given by $x \supset y := x \imp (x \land y)$ that will be used
for our twist representation. 
Unlike the case of $\{\to, \nnot\}$, this fragment determines a sub-classical logic coinciding with the deductive system  known in the literature 
as \emph{Soboci\'nski's logic}, which is defined by the matrix $\la \Al[S]_3, \{\uv, \tv \}\ra$, where 
$\Al[S]_3$ denotes the $\{\imp, \nnot\}$-reduct of the three-element OL-algebra \cite{sobocinski,Pa72}. (Note that $\Al[S]_3$ is not itself an OL-algebra because,
for instance, its operations are not able to define the constant $\uv$.)
%
Soboci\'nski originally proposed
a system 
that would be able to avoid the paradoxes of material implication, its main interest 
consisting in the fact that no formula in which an atomic proposition appears once is satisfied by the above matrix~\cite{sobocinski}. This intends to cover famous examples such as `if $\varphi$ then (if $\psi$ then $\varphi$)'.
\\

\begin{figure}[tbh]
    \centering
    \begin{tabular}{@{}c|ccc@{}}
         \toprule
          $\imp$ &  \bv & \uv & \tv \\
         \midrule
          \bv & \tv & \tv & \tv \\
          \uv & \bv &  \uv &\tv\\
          \tv &  \bv & \bv & \tv \\
         \bottomrule
     \end{tabular}
     \quad
     \begin{tabular}{@{}c|ccc@{}}
     \toprule
         $\supset$ &  \bv & \uv & \tv \\
         \midrule
          \bv & \tv & \tv & \tv \\
          \uv &\bv & \uv &\tv\\
          \tv & \bv & \tv &\tv \\
        \bottomrule
    \end{tabular}
    \quad
     \begin{tabular}{@{}c|ccc@{}}
     \toprule
         $\otimes$ &  \bv & \uv & \tv \\
         \midrule
          \bv & \bv & \bv & \bv \\
          \uv &\bv & \uv &\tv\\
          \tv & \bv & \bv &\tv \\
        \bottomrule
    \end{tabular}
\caption{Truth tables of Soboci\'nski connectives.}
\label{fig:truth-tables-deriv}
\end{figure}
As noted in~\cite{Pa72},
the valid formulas of Soboci\'nski's logic (but not the whole consequence relation) also coincide with those of the corresponding fragment
of the relevance logic known as R-mingle.
%
Let us call a \emph{Soboci\'nski algebra} any member of the quasi-variety 
$Q(\Al[S]_3)$\footnote{An abstract presentation of $Q(\Al[S]_3)$, though not explicitly given, can easily be obtained 
from the results in~\cite{BlRa08}; another presentation is implicit in Proposition~\ref{prop:sobofirsancheck} below.}. 
Blok and Raftery~\cite{BlRa08} 
observe that $Q(\Al[S]_3)$ is not a variety\footnote{The variety $V(\Al[S]_3)$ is itself an interesting class: on the one hand, its associated $\tau$-assertional logic (with $\tau (x) := x \imp x$) is not protoalgebraic; on the other, it is not know whether $V(\Al[S]_3)$ can be finitely axiomatized 
(see~\cite{BlRa08} and~\cite[Ex. 9]{Ra06}).}.
Indeed, the congruence $\theta$ that identifies $\bv$ with $\tv$ 
determines a two-element quotient algebra $\Al[S]_3 / \theta$ that does not
satisfy 
the 
symmetric version of our~\eqref{eq:deford} introduced below; 
we also note that $\Al[S]_3 / \theta$ instead satisfies the interesting equation
$x = \nnot x$. We shall return on this observation below, once we are able to take advantage of 
the twist representation. 
\\

Our next goal is precisely to show that all Soboci\'nski algebras can be viewed as twist-algebras in the sense
of the next definition. 

\begin{definition}
\label{d:sobtwist}
Let
$ \Al[B]  = \la B ; 
\wedge , 
\to ,  
1
\ra$
be a generalized Boolean algebra.
Define the algebra 
$ 
\Al[B]^{\bowtie}
:= \la B^{\bowtie}; 
\imp, \nnot
\ra$ 
with universe 
$B^{\bowtie}  : =  \{ \la a_1 , a_2  \ra \in B \times B : a_1 \to a_2 = a_2 
\}
$
and 
operations given,
for all $\la a_1 ,a_2 \ra , \la b_1, b_2 \ra \in B \times B$, by:
\begin{enumerate}[(i)]

\item $\nnot \la x_1 , x_2 \ra := \la x_2,   x_1  \ra$, 
\item $\la x_1  , x_2 \ra \imp \la y_1 , y_2 \ra  := \la  (x_1 \to y_1) \land (y_2 \to x_2), x_1  \land y_2  \ra$. 
\end{enumerate}
%
A \emph{
\sobocinski twist-algebra over $\alg{B}$}
is any subalgebra $\Al[A] \leqslant \Al[B]^{\bowtie} 
$
satisfying
$\pi_1[A] = B$. 
\\

\end{definition}

It is obvious from the two constructions that 
every Soboci\'nski twist-algebra embeds into an OL-twist algebra: thus we have that 
Soboci\'nski twist-algebras are precisely  the $\{\imp, \nnot\}$-subreducts of OL-algebras.
This, in turn, entails that each Soboci\'nski twist-algebra $\Al[A] \leqslant \Al[B]^{\bowtie}$ 
belongs to $Q(\Al[S]_3)$.
It is also clear that 
$ \Al[S]_3$ may be viewed as the (unique)
Soboci\'nski twist-algebra over the two-element Boolean algebra,
with $\tv := \la 1, 0 \ra $, $\uv :=  \la 1, 1 \ra $ and $\bv := \la 0, 1 \ra$.
\\

On each twist-algebra, the OL-conjunction operation $\land$ can be defined by the term:
$
x \land y := \nnot (x \imp \nnot y)
$
and, vice versa, one could define
$x \imp y := \nnot (x \land \nnot y).$
On pairs, the preceding definition gives us:
$$
\la x_1  , x_2 \ra \land \la y_1 , y_2 \ra  = \la x_1  \land  y_1 , ( x_1 \to  y_2) \land  ( y_1 \to  x_2) \ra.
$$
As on OL-algebras, $\land$ is a semilattice operation.
However, on each
Soboci\'nski twist-algebra $\Al[A] \leqslant \Al[B]^{\bowtie}$, a different canonical order $\leq $ (not coinciding with the one induced by $\land$)
is given by:
\begin{equation}
\label{eq:deford}
\tag{$\dagger$}
x \leq y \qquad \text{ if and only if } \qquad x \imp y = (x \imp y) \imp (x \imp y),
\end{equation}
while on pairs we have that
$
\la a_1  , a_2 \ra \leq  \la b_1 , b_2 \ra \text{ if and only if } a_1 \leq b_1 \ \text{ and } \ b_2 \leq a_2.
$
%
%
\\

Now, let us further define the following terms (cf.~\cite[p.~647]{BlRa08}):
$$
\abs{x} : = x \imp x \qquad 
\text{ and }\qquad  
x \supset y : = (x \imp \abs{y} ) \imp (x \imp y).
$$
We shall be especially interested in the implication $ \supset$, which can equivalently be defined
as follows:
$$x \supset y : = x \imp (x \land y).
$$

%
\begin{remark}
    \label{remark:cioresob}
    The connective $\supset$ coincides with the primitive implication
of the three-valued logic of formal inconsistency known as \emph{Ciore} or \emph{LFI2}
(see, e.g., \cite{Co21,CaMa00}); we shall say more on this
in a while.
\\

\end{remark}

\begin{fact}
\label{fact:imppairs}
Given a Soboci\'nski twist-algebra $\Al[A] \leq \Al[B]^{\bowtie}$ and $\la a_1, a_2 \ra, \la b_1, b_2 \ra \in A$, we have:
\begin{align*}
\la a_1  , a_2 \ra \supset \la b_1 , b_2 \ra  &  = \la  a_1  \to  b_1  ,   a_1 \land  b_2 \land (b_1 \to a_2) \ra.
\end{align*}

\begin{proof}
Since 
    \begin{align*}
 \la a_1, a_2 \ra \supset \la b_1, b_2 \ra 
 & = \la (a_1 \to (a_1  \land  b_1)) \land ((( a_1 \to  b_2) \land  ( b_1 \to  a_2)) \to a_2) , a_1 \land ( a_1 \to  b_2) \land  ( b_1 \to  a_2) \ra \\
 & = \la (a_1 \to  b_1) \land ((( a_1 \to  b_2) \land  ( b_1 \to  a_2)) \to a_2) , a_1 \land   b_2 \land  ( b_1 \to  a_2) \ra,
\end{align*}
it suffices to show that 
$a_1 \to  b_1 \leq (( a_1 \to  b_2) \land  ( b_1 \to  a_2)) \to a_2$. By residuation, the latter is equivalent to
$(a_1 \to  b_1) \land ( a_1 \to  b_2) \land ( b_1 \to  a_2) \leq  a_2$.
To conclude the proof, it suffices to recall the assumption $a_1 \to a_2 = a_2$ and to note that 
$(a_1 \to  b_1) \land ( b_1 \to  a_2) \leq  a_1 \to a_2 = a_2$.
\end{proof}
    
\end{fact}

The component-wise definition of $\supset$ suggests that, within OL, the formula $\phi \supset \psi$ is inter-derivable with primitive implication 
(as well as with the term-definable weak Nelson implication),
though their negations need not be. Either of these connectives in fact provides a Deduction Detachment Theorem
for OL  (but for Soboci\'nski's logic one must rely only on $\supset$).
For our purposes, the main importance of 
the operation $\supset$ is
that it allows us to recover, on the first component, 
precisely the Boolean 
implication. This is the key to the twist representation. 
Writing $x \pr y$ whenever $x \supset y = \abs{x \supset y}$, 
it is easy to check that: 
$$
\la a_1  , a_2 \ra \pr  \la b_1 , b_2 \ra \qquad \text{ iff } \qquad a_1 \leq b_1.
$$
By~\eqref{eq:deford}, we have that
$x \pr y$ \ iff \ $x \leq x \land y$.
It is also clear that $x \leq y$ holds if and only if both $x \pr y $ and $\nnot y \pr \nnot x$ hold. 
\\

The  properties listed below are clearly satisfied by every Soboci\'nski twist-algebra.
As for general Soboci\'nski algebras,
 these being quasi-equational properties, they can be verified directly on 
 $\Al[S]_3$:

\begin{proposition}
\label{prop:sobofirsancheck}
The following properties hold on every  Soboci\'nski algebra:
\begin{enumerate}[(i)]




\item $\pr$ is a preorder (i.e., reflexive and transitive),


\item $x \pr y$ implies $x \land z \pr y \land z$, $y \supset z \pr x \supset z$ and $z \supset x \pr z \supset y$,
 \item  In consequence, $\eq$ is a congruence  on $\la A; \land, \supset \ra$,

\item The quotient $\la A; \land, \supset \ra/ \eq $ is a generalized Boolean algebra,

\item $x \leq y$ (i.e.~$x \imp y = \abs{x \imp y}$) \ if and only if \ [$x \pr y$ and $\nnot y \pr \nnot x]$,

\item $ x \imp y \equiv (x \supset y) \land (\nnot y \supset \nnot x)$,

\item $\nnot x \eq x  \supset \nnot x    $.
\end{enumerate}
\end{proposition}

%
%
%
%

\begin{theorem}[Soboci\'nski twist representation]
\label{d:sobtwistrep}
Every Soboci\'nski algebra $\A$ 
embeds into the twist-algebra
$[\la A; \supset \ra / \eq]^{\bowtie}$ constructed according to Definition~\ref{d:sobtwist}.
\end{theorem}

\begin{proof}
Proposition~\ref{prop:sobofirsancheck} ensures that we have the correct ingredients to build a twist-algebra. 
The embedding $\iota \colon A \to A/\!\!\eq \times A/ \eq $ is defined, as expected, by
$\iota(a) := \la [a], [\nnot a] \ra $, and it is
is injective by  item (v) of Proposition~\ref{prop:sobofirsancheck}.
Let us check that $\iota$ preserves the algebraic operations. 
The case of the negation 
immediately follows from the involutive law. As for the implication, we have:
\begin{align*}
\iota (a \imp b) 
& = \la [a \imp b], [ \nnot (a \imp b)] \ra \\
& = \la [a \imp b], [ \nnot \nnot (a \land \nnot b)] \ra \\
%
%
& =  \la  [(a \supset b)  \land (\nnot b \supset \nnot a)], [a \land \nnot b] \ra & \text{by~ Prop.~\ref{prop:sobofirsancheck} (vi)}\\
& =  \la  ([a] \to [b ]) \land ([\nnot b] \to [\nnot a]), [a] \land [\nnot b] \ra \\
& =  \la  [a], [\nnot a ] \ra \imp  \la  [b], [\nnot b] \ra \\
& = \iota (a) \imp \iota (b).
\end{align*}
Finally, observe that
item (vii) of Proposition~\ref{prop:sobofirsancheck}
ensures that,
for every element $\la [a], [\nnot a] \ra \in \iota[A]$, we have
$[a ] \to [\nnot a ] = [a \supset \nnot a]= [ \nnot a ] $, as required by Definition~\ref{d:sobtwist}.
\end{proof}

Joining the previous results, we have that Soboci\'nski twist-algebras coincide, up to isomorphism, with the members of 
$Q(\Al[S]_3)$. This suggests a straightforward way to embed an arbitrary algebra in $Q(\Al[S]_3)$
into an OL(-twist)-algebra;
note also that Proposition~\ref{prop:sobofirsancheck} contains a set of quasi-equations
that can be easily adapted to provide a complete presentation for $Q(\Al[S]_3)$. 
\\

 
    The following result generalizes and clarifies what we observed earlier about the algebra $\Al[S_3]$
    and the congruence identifying $\bv$ and $\tv$, and gives some information
    on those algebras in $V(\Al[S_3])$ that are not in $Q(\Al[S_3])$.
Given a Soboci\'nski algebra $\A$, let $\theta$ denote the congruence 
generated by the set $\{ \la a, \abs{a }\ra : a \in A \} $. 

    \begin{lemma}
    \label{l:congabs}
Let $\A$ be  Soboci\'nski algebra. 
Then $\theta = \{ \la a, b \ra : \abs{a} = \abs{b} \}$ and
 $\Al / \theta \notin Q(\Al[S_3])$ for  each   non-trivial $\Al / \theta$.
\end{lemma}

\begin{proof}
Let $\theta$ be  the congruence 
generated by the set $\{ \la a, \abs{a }\ra : a \in A \} $,
and let
$R : = \{ \la a, b \ra : \abs{a} = \abs{b} \}$.
If
$\la a,b \ra \in R$, 
then
$\abs{a} \theta \abs{b}$.
We thus have
$a \theta \abs{a} \theta \abs{b} \theta b$,
showing that 
$
R \subseteq \theta$.
To show the converse inclusion, it suffices to verify that 
$\{ \la a, b \ra : \abs{a} = \abs{b} \}$ is a congruence
containing every pair of type $\la a, \abs{a} \ra$.
The latter claim follows from the observation that the equation
$\abs{x} = \abs{x} \imp \abs{x}$ is valid on $\Al[S_3]$.
To show that $R$ is a congruence,   assume
 $\Al[A] \leq \Al[B]^{\bowtie}$ is a twist-algebra. 
Note that 
$\abs{\la a_1, a_2 \ra} = \abs{\la b_1, b_2 \ra} $ if and only if $a_1 \land a_2 = b_1 \land b_2$.
This immediately implies that $R$ 
is  an equivalence relation compatible with the negation; note, furthermore, that 
$\la \nnot a,   a \ra \in R$ 
for all $a \in A$.
As to the implication operation, letting
$\la c_1, c_2 \ra \in A$, we have:
$\la a_1, a_2 \ra \imp \la c_1, c_2 \ra = \la (a_1 \to c_1) \land (c_2 \to a_2), a_1 \land c_2 \ra $, so that 
$(a_1 \to c_1) \land (c_2 \to a_2) \land  a_1 \land c_2 \ra  = a_1 \land c_1 \land c_2 \land a_2 = b_1 \land c_1 \land c_2 \land b_2 $.
Hence, $\la a_1, a_2 \ra \imp \la c_1, c_2 \ra R \la b_1, b_2 \ra \imp \la c_1, c_2 \ra$. A similar computation 
shows that $\la c_1, c_2 \ra \imp \la a_1, a_2 \ra R \la c_1, c_2 \ra \imp \la b_1, b_2 \ra$.
Thus $R = \theta$.
Finally, since
each quotient $\Al / \theta $ will satisfy 
$x = \abs{x}$,
assuming $\Al / \theta \in Q(\Al[S_3])$, we would have
$a \imp b = \abs{a \imp b}$ and 
$b \imp a = \abs{b \imp a}$ for all $a, b \in \Al / \theta$, entailing  $a =b$, by \eqref{eq:deford}. Hence, $\Al / \theta$ must be trivial. 
\end{proof}

Another way to obtain quotient algebras that are 
not in $ Q(\Al[S_3])$ results from the following construction,
which is inspired by the examples shown in~\cite[Sec.~4.1]{BlRa04}.
Given a  Soboci\'nski (twist-)algebra $\Al \in Q(\Al[S_3])$ 
and a subset $B \subseteq A$,
define the set
$B^* := \{ a \in A : b \leq \abs{a} \text{ for some } b \in B\}$,
where
the order $\leq$ is given by
$x \leq y$ iff $x \imp y = \abs{x \imp y}$.
On a twist-algebra, 
we have that $\la a_1, a_2 \ra \in B^*$
if and only if 
there is
$\la b_1, b_2 \ra \in B$
such that 
$a_1 \land a_2 \leq b_2$.
Thus, for singletons $B = \{ b \}$ and
$C = \{ c \}$, we
have $B^* \subseteq C^*$
if and only if 
$\nnot b \pr \nnot c$.
In particular, if
$  \la b_1, 1 \ra \in B $, 
we simply have
$B^* =A$, while for
$B = \{ \la b_1, 0 \ra \} $,
the set 
$B^*$
consists precisely of the Boolean elements of $\A$,
i.e.~those $\la a_1, a_2 \ra$ such that 
$\neg a_1 = a_2$.
The following lemma shows how to use 
$B^*$
to define a congruence $\theta_B$ such that $\A / \theta_B \notin Q(\Al[S_3])$.

\begin{lemma}
\label{l:nonboolel}
Let $\Al \in Q(\Al[S_3])$ and $B \subseteq A$, with $B^*$ 
defined as before.
\begin{enumerate}[(i)]
    \item $b \imp a, a \imp b \in B^*$ all $b \in B^*$ and $ a \in A$.
    \item In consequence, $B^*$
    is the universe of a 
    subalgebra of $\A$.
    \item $\theta_B : = Id_A \cup (B^* \times B^*)$ is a congruence on $\A$.
    \item 
    $
    \Al / \theta_B \notin Q(\Al[S_3])$ for each non-trivial $\Al / \theta_B$.
\end{enumerate}
\end{lemma}
\begin{proof}
(i). Assume $\Al[A] \leq \Al[B]^{\bowtie}$ is a twist-algebra.
As observed above, $\la b_1 , b_2 \ra \in B^*$ if and only if
$b_1 \land b_2 \leq c_2$ for some $\la c_1, c_2 \ra\in B$.
%
%
Then, computing
$\la a_1  , a_2 \ra \imp \la b_1 , b_2 \ra  = \la  (a_1 \to b_1) \land (b_2 \to a_2), a_1  \land  b_2  \ra$,
it suffices to note that
$(a_1 \to b_1) \land (b_2 \to a_2) \land a_1  \land  b_2 = a_1 \land a_2 \land b_1 \land b_2 \leq c_2$.
A similar argument establishes the case of 
the implication
$\la b_1 , b_2 \ra \imp \la a_1  , a_2 \ra$.

Clearly, the preceding argument also establishes item (ii); note that $B \subseteq B^*$
but $B^*$ may be larger, in general, than the subalgebra generated by $B$.

(iii). Suffice it to note that, if $\la a, b \ra \in \theta_B$ are such that $a \neq b$, then $a,b \in B^*$.
Then, by item (i), we have $a \imp c, b \imp c, c \imp a, c \imp b \in  B^*$, for all $c \in A$.
So all these elements are related by $\theta^B$, as required. 

(iv). By construction, all elements
in $B^* / \theta_B$ are designated. We
may assume that $\A / \theta_B$ has at least an element that is not designated:
otherwise, we would conclude $\A / \theta_B \notin Q(\Al[S_3])$
already from Lemma~\ref{l:congabs}. Let then $a/ \theta^B$ be such a non-designated element (for some $a \in A$).
Recall that every algebra in $Q(\Al[S_3])$ satisfies the quasi-equation corresponding to \emph{modus ponens}, that is, $x = \abs{x}$ and $x \imp y = \abs{x \imp y}$ entail $y = \abs{y}$. 
Now, this
quasi-equation fails in $\Al / \theta_B$.
Indeed, choosing any $b \in  B^*$,  we have that $b/ \theta_B \imp a/ \theta_B = b/ \theta_B$ is designated while $a/ \theta_B$ is not. 
\end{proof}

\begin{example}
Let $\Al[B]_4$ be the four-element Boolean algebra with universe 
$B_4 = \{ 0, 1, a, \neg a\}$. The twist-algebra
$\Al[B]_4^{\bowtie}$ (isomorphic to the direct power
$\Al[S]_3 \times \Al[S]_3$) has for universe
the nine-element set: 
$$
\{ \la 0, 1 \ra, \la 1, 0 \ra, \la 1, 1 \ra, \la a, \neg a \ra, \la \neg a, a \ra, 
\la 1, a \ra, \la a, 1 \ra, \la 1, \neg a \ra, \la \neg a, 1 \ra  \}.
$$
Considering
$B = \la 1, 0 \ra$, we have
$B^* =
\{ \la 0, 1 \ra, \la 1, 0 \ra, \la a, \neg a \ra, \la \neg a, a \ra \}
$.
Thus, the quotient algebra $\Al[B]_4^{\bowtie} / \theta_B$
obtained via the congruence considered in Lemma~\ref{l:nonboolel}
has for universe the six-element set: 
$$
\{ [\la 1, 0 \ra], [\la 1, 1 \ra],
[\la 1, a \ra], [\la a, 1 \ra], [\la 1, \neg a \ra], [\la \neg a, 1 \ra]  \}.
$$
The designated elements on this algebra are
$
\{ [\la 1, 0 \ra], [\la 1, 1 \ra],
[\la 1, a \ra], [\la 1, \neg a \ra]  \},
$
and
the failure of the quasi-equation corresponding to \emph{modus ponens}
is witnessed, for instance, by the equality 
$[\la 1, 0 \ra] \imp [\la a, 1 \ra] = [\la 1, 0 \ra]$.
We note that, among the subalgebras of $\Al[B]_4^{\bowtie} / \theta_B$,
the four-element one 
with universe
$
\{ [\la 1, 0 \ra], [\la 1, 1 \ra],
[\la 1, a \ra], [\la a, 1 \ra] \}
$
is isomorphic precisely to the algebra $\Al[Z_3^*]$ considered in~\cite[Lemma~27]{BlRa04}.
\end{example}

We end our discussion by briefly returning on the comparison between Soboci\'nski's logic and Ciore (Remark \ref{remark:cioresob}) in the light of our twist representation. Ciore is defined from a three-element algebra whose  operations 
(conjunction, disjunction, implication and negation)
can be obtained from those of $\Al[S]_3$
as follows:
\begin{itemize}
\item the conjunction and the negation coincide with those of $\Al[S]_3$;
\item the implication is the above-defined operation $\supset$;
\item the disjunction corresponds to the following $\Al[S]_3$-term\footnote{This is the term denoted by
$\lambda (x,y)$ in~\cite[p.~79]{BlRa04},
 also independently rediscovered in~\cite{Greati23} and employed for the purpose of axiomatizing
 the $\{\imp, \nnot\}$-fragment of OL. 
}: $ ((x \imp y) \land (y \imp x)) \imp y$.
\end{itemize}
Actually, one can verify that the conjunction (and therefore also the disjunction) can be defined 
in the $\{ \supset, \nnot \}$-fragment, 
as follows
(see Figure~\ref{fig:truth-tables-deriv}):
\begin{align*}
x \otimes y & := \nnot (x \supset \nnot y),\\
x \land  y & := \nnot ((x \supset \nnot y ) \otimes ( y \supset \nnot x) ).
\end{align*}
(This also entails that the $\{ \supset, \nnot \}$-fragment of Soboci\'nski's logic already expresses
all the Soboci\'nski connectives).
In addition to the above, Ciore also has a unary ``consistency'' connective $\circ$ given by 
$\circ \uv = \bv$ and $\circ \bv = \circ \tv = \tv$. This operator is not definable in the language of 
$\Al[S]_3$ unless we expand it with the constant $\bv$: then one can verify that
$\circ x  = (x \land \nnot x) \imp \bv$. Since the designated elements in Ciore are the same as in
Soboci\'nski's logic, we conclude that the former is precisely the conservative expansion of the latter
obtained by adding the constant $\bv$.

\subsection*{The $\{\imp, \bv\}$-fragment}

As a technical curiosity, let us see that, although the twist construction usually hinges on the 
involutive negation, 
it is also possible to adapt it to give a representation for some fragments that lack it.
We will just sketch the main ideas of the construction for
the case of 
the $\{\imp, \bv\}$-fragment of Ciore (or Soboci\'nski--OL, if expanded with the constant $\bv$). 

Consider   a twist-algebra 
$\Al[A] \leqslant \Al[B]^{\bowtie}$ (in the sense of Definition~\ref{d:sobtwist})
such that $\bv = \la 0, 1 \ra \in A$. 
Note that, in this case, we also have $ \bv \imp \bv= \tv = \la 1, 0 \ra \in A$. Let us introduce the following terms:
$$
   \Box  x 
   : = \tv \imp x, \qquad
         \Diamond  x 
         : = x \imp \bv, \qquad
   x \odot y 
   := \Diamond \Box (x \imp \Diamond y).
$$
Let $\la a_1, a_2 \ra, \la b_1, b_2 \ra \in A$. Note that
$\neg a_1 := a_1 \to 0 \leq a_1 \to a_2 = a_2$, and, symmetrically,
$\neg a_2 \leq a_1$.
Thus, computing component-wise, we have:
\begin{align*}
   \Box  \la a_1, a_2 \ra &  = \la (1 \to a_1) \land (a_2 \to 0), 1 \land a_2 \ra  = \la \neg a_2, a_2 \ra,  \\
        \Diamond \la a_1, a_2 \ra &  =  \la (a_1 \to 0) \land (1 \to a_2), a_1 \land 1 \ra  = \la \neg a_1, a_1 \ra,  \\
\Diamond \Box \la a_1, a_2 \ra & = \la a_2, \neg a_2 \ra,  \\
\Diamond \Diamond  \la a_1, a_2 \ra & = \la  a_1, \neg  a_1 \ra, \\
   \la a_1, a_2 \ra \odot \la b_1, b_2 \ra & =  \Diamond \Box \la (a_1 \to \neg b_1) \land (b_1 \to a_2), a_1 \land b_1 \ra = 
   \la a_1 \land b_1, \neg (a_1 \land b_1) \ra.
\end{align*}
The operators $\Box$ and $\Diamond$ will jointly allow us to somehow
replicate the required features of
the usual twist-algebra negation.
Note that in this fragment
 Ciore's consistency operator $\circ$ is also definable (e.g.~by $\circ x := \Box 
 (x \imp  x)
 $).
The component-wise definitions entail that every Soboci\'nski (twist-)algebra with $\bv$ satisfies the following properties
(we use the abbreviation 
$\eq$ 
as before):
\begin{enumerate}[(i)]
\item $x \eq y $  \ iff \ $\Diamond x = \Diamond y$ \ iff \ $\Diamond \Diamond x = \Diamond \Diamond y$.
\item $\nnot x \eq \nnot y $  \ iff \ $\Box x = \Box y$ \ iff \ $\Diamond \Box x = \Diamond \Box y$. 
\item In consequence, $x = y $ \ iff \ [$\Diamond \Diamond x = \Diamond \Diamond y$ and $\Diamond \Box x = \Diamond \Box y$].
\item $\Diamond \Diamond (x \imp y) = 
(\Diamond \Diamond a \imp \Diamond \Diamond b) \odot 
(\Diamond \Box b \imp \Diamond \Box a ) $.
\item $ \Diamond \Box (x \imp y) = \Diamond \Diamond x \odot \Diamond \Box y   $.
\end{enumerate}

Let us further notice that the forward image $A_{\Diamond} := \{ \Diamond \Diamond a : a \in A \}$ is closed under 
the operations $\{\odot, \imp, \bv, \tv \}$, so that 
$\Al_{\Diamond} = \la A_{\Diamond} , \odot, \imp, \bv, \tv \ra $ is a Boolean algebra.
The idea is thus to 
embed an algebra over the $\{\imp, \bv\}$-language satisfying the above properties
into a twist-algebra (built as per Definition~\ref{d:sobtwist}, including the constant
$\bv = \la 0, 1 \ra $ 
but disregarding the negation)
over $\Al_{\Diamond} $.
For this, one can define the map $\iota (a) := \la \Diamond \Diamond a,   \Diamond \Diamond \Diamond \Box a \ra $ 
(or simply $\iota (a) := \la \Diamond \Diamond a,    \Diamond \Box a \ra $,
for
$\Diamond  x =\Diamond \Diamond \Diamond  x $), 
using the above properties to verify that (the negation and) the implication are preserved
as follows:
\begin{align*}
\iota (a \imp b) 
& = \la \Diamond \Diamond (a \imp b),  \Diamond \Box (a \imp b) \ra \\
%
%
%
& =  \la  (\Diamond \Diamond a \imp \Diamond \Diamond b) \odot 
(\Diamond \Box b \imp \Diamond \Box a ),  \Diamond \Diamond a \odot \Diamond \Box b \ra \\
& =  \la  \Diamond \Diamond a, \Diamond \Box a  \ra \imp  \la  \Diamond \Diamond b, \Diamond \Box b \ra \\
& = \iota (a) \imp \iota (b).
\end{align*}

\subsection{De Finetti's logic}
\label{ss:df}

De Finetti's approach for dealing with indicative conditionals \cite{deFin} built on previous considerations by Reichenbach (see, e.g., \cite{egre20211,wollic} for more details). The guiding thought experiment here relies on probabilistic intuitions. Roughly, De Finetti sees propositions as events and defines the event corresponding to $\varphi \to \psi$ as the \emph{tri-event} `$\psi$ given $\varphi$' or, informally, as a conditional bet on the truth-value of $\varphi$ given $\psi$: it is won in case both $\varphi$ and $\psi$ are fulfilled, lost when $\varphi$ is fulfilled but $\psi$ is not and simply called off when $\varphi$ is not fulfilled. This is consistent with the behavior of the connectives that we have already presented before as an example: 
$$\varphi \to_\mathsf{DF} \psi = \tv \text{ iff } \varphi = \tv \text{ and } \psi = \tv, \quad \varphi \to_\mathsf{DF} \psi = \bv \text{ iff } \varphi = \tv \text{ and } \psi = \bv,$$
and, additionally,
$$\varphi \land_\mathsf{K} \psi = \tv \text{ iff } \varphi = \tv \text{ and } \psi = \tv, \quad \varphi \land_\mathsf{K} \psi = \bv \text{ iff } \varphi = \bv \text{ or } \psi = \bv.$$
Some interesting features of DF are already known. For instance, DF is truth-equational, but fails to be self-extensional and protoalgebraic (thus algebraizable) \cite[Props. 1,2]{wollic}. In fact, DF is known to be an expansion of Priest's logic of paradox \cite{priest} with constant $\uv$ \cite{wollic}. On the other hand, the algebraic models of DF can be seen as centered Kleene lattices (Definition \ref{kleene}). Specializing a construction from \cite{um20} we can define the following:
\begin{definition}
\label{DFtwist}
     Let $\*L := (L; \land, \lor, 1)$ be an upper-bounded distributive lattice. The \emph{full DF-twist-algebra over $\mathbf{L}$} consists in an algebra $\mathbf{L}^{\bowtie}$ with universe $\{ (a_1 , a_2) \in L \times L \, | \,  a_1 \lor a_2 = 1 \}$ and operations
\begin{itemize}
    \setlength \itemsep{0pt}   
    \item[i.] $(x_1, x_2) \land (y_1, y_2) := (x_1 \land y_1, x_2 \lor y_2)$,
    \item[ii.] $\neg(x_1, x_2) := (x_2, x_1)$, 
    \item[iii.] $\uv :=  (1, 1)$.
\end{itemize}
We can then define $x \lor y : = \nnot (\nnot x \land \nnot y)$ and 
$x \to y := (\uv \land \nnot x) \lor (x \land y).$ A \emph{DF-twist-algebra over $\*L$} is any subalgebra $\*A\leqslant {\*L}^{\bowtie}$ such that $\pi_1[A] = L$. 
\end{definition}
Then, every DF-twist-algebra is a centered Kleene lattice. Let us note that we impose the requirement $a_1 \lor a_2 =1 $ instead of  $a_1 \land a_2 =0$, as it is done in the Kallman construction for centered Kleene algebras \cite{Cig86}. This is connected to the fact that we take the center to be $(1, 1)$ instead of $(0, 0)$. Later on we will consider \textit{both} pairs. We can make the following observations: 

\begin{remark}
\label{toDF}
One can check that, on DF-twist-algebras, 
$$(x_1, x_2) \to_\mathsf{DF} (y_1, y_2) = (x_2 \lor (x_1 \land y_1), x_2 \lor y_2).$$
\end{remark}

\begin{fact}
\label{notfullDF}
    A DF-twist-algebra over $\*L$ is not necessarily full.

\begin{proof}
    Indeed, take $\*L$ to be the distributive lattice given by $1 < a, b < 0$. Then, $\*L^{\bowtie}$ will contain the elements $(1, 1), (1, a), (1, b), ( 1, 0), ( a, b)$ together with their symmetrizations. Now, consider the subset $A$ obtained by leaving aside $( a, b)$ and $( b, a)$. Then, note that $A$ is closed under negation, contains $\uv$ and $\pi_1(A) = L$. It remains to see that it is also closed under conjunction in order to conclude that it is the desired subuniverse of $\*L^{\bowtie}$. But this is easy to check by assuming that $( a, b)$ is expressible in terms of the remaining elements and arriving at a contradiction.
\end{proof}
\end{fact}
Now, given a centered Kleene algebra $\*A := (A; \land, \lor, \neg, \uv)$, define the map $\Diamond: A \to A: x \mapsto x \land \uv$. Then, we define the algebra $\Diamond(\*A) := (\Diamond[A]; \land, \lor, \uv)$, in which the operations are defined as restrictions of the ones from $\*A$. 
\begin{lemma}
    \label{DFlemma}
    Given a centered Kleene algebra $\*A$, the algebra $\Diamond(\*A)$ is an upper-bounded distributive lattice. 
\end{lemma}
\begin{proof}
   We only need to note that $\Diamond \Diamond x = \Diamond x$ and that $\Diamond$ preserves both $\land$ and  $\lor$ (by distributivity).  
\end{proof}

\begin{theorem}[DF-twist representation]
\label{DFtwistrep}   
Every centered Kleene algebra $\*A$ is isomorphic to a DF-twist-algebra over $\Diamond(\*A)$ via the map  $\iota \colon \*A \to \Diamond (\*A) \times \Diamond (\*A)$ given by $ a \mapsto (\Diamond a, \Diamond \nnot a)$.
\end{theorem}
\begin{proof}
To verify that the map $\iota$ is injective, observe that,
for all $a,b \in A$, we have
$a \lor \uv = b \lor \uv$ if and only if 
$\Diamond \nnot a  = \Diamond \nnot b $.
Using this observation,  assuming $\iota(a) = \iota (b)$, we have:
\begin{align*}
    a 
    & = a \land (a \lor \uv) & \text{by absorption} \\
    & = a \land (b \lor \uv) &  \Diamond \nnot a  = \Diamond \nnot b 
    \\
    & = (a \land b) \lor (a \land \uv) & \text{by distributivity} \\
    & = (a \land b) \lor (b \land \uv) & \Diamond a = \Diamond b 
    \\
    & = (b \land a) \lor (b \land \uv) & \text{by commutativity} \\
    & = b \land (a \lor \uv) & \text{by distributivity} \\
    & = b \land (b \lor \uv) & \Diamond \nnot a  = \Diamond \nnot b 
    \\
    & = b & \text{by absorption}.
\end{align*}
The definition of DF-twist-algebras also requires
$\Diamond a \lor \Diamond \nnot a = \uv$ for all $a \in A$.
To see this, recall we can instantiate the Kleene equation as
$\uv = \uv \land \nnot \uv \leq a \lor \nnot a $.
Then, we have:
\begin{align*}
 \Diamond a \lor \Diamond \nnot a 
  & =   ( a \land \uv) \lor ( \nnot a \land \uv ) \\   
  & =    ( a \lor \nnot a) \land \uv  & \text{by distributivity}\\ 
  & = \uv.
\end{align*}
Let us check that $\iota$ preserves the algebraic operations.
The case of the constant is trivial, and negation is 
straightforward:
for all $a \in A$, we have
$$\iota (\nnot a)  
= (\Diamond \nnot a,  \Diamond \nnot \nnot a)
= \nnot (\Diamond a, \Diamond \nnot a)
= \nnot \iota (a).$$
Regarding the meet, we have:
\begin{align*}
\iota (a \land b) 
& = ( \Diamond (a \land b), \Diamond \nnot (a \land b) )\\
& = (  a \land b \land \uv,  \nnot (a \land b) \land \uv )\\
& = (  a \land b \land \uv,  (\nnot a \lor \nnot b) \land \uv )& \text{by~De Morgan's law} \\
& =  (  a \land \uv \land   b \land \uv,   (\nnot a \land \uv ) \lor  ( \nnot b \land \uv) )
& \text{by~distributivity} \\
& =  ( \Diamond a \land  \Diamond b,  \Diamond \nnot a \lor  \Diamond \nnot b  )\\
& =  ( \Diamond a, \Diamond \nnot a  )\land  ( \Diamond b, \Diamond \nnot b  )\\
& = \iota (a) \land \iota (b).
\end{align*}
\end{proof}

\subsection{The logics of Farrell and Cantwell}

The case of F and CN is interesting both in a logico-mathematical and a philosophical sense. Not only these two logical systems were introduced within a difference of more than twenty years from each other but their original motivations differed completely. Despite this, it turns out that F and CN are definitionally equivalent, i.e. from the formal standpoint, both logics are essentially the same. This holds in virtue of the following equalities \cite[Thm. 9]{wollic}:
\begin{align*}
    x \OLimp y \ &  : = \, \OLneg ((y \Fimp x) \Fimp \OLneg y) \Fimp ((x \uor y) \Fimp y), \\
    x \Fimp y  \  & : =   \, x \OLimp (x \uand y).
\end{align*}
However, later we will see that the four-valued extensions of F and CN need not to behave in the same way. This is why we will present both logics in a differentiated way below, even if we will make use of the previous fact.

\subsection*{Farrell}

Farrell's motivations \cite{Fa79,Fa86} depart from Strawson's theory of presupposition and the study of a satisfying theory of quantification. Even though Farrell provides a technical justification for the table of implication, he intends to see $\uv$ as the truth-value `to be ignored' \cite{Fa86}, similarly as in programming languages in which the evaluation of a formula with false antecedent is omitted. Finally, Farrell wished to deploy his system in order to elucidate philosophical problems such as Hempel's paradox \cite{Fa86}. The conditional of this logic has a mixed behavior between those of De Finetti and Cooper: 
$$\varphi \to_\mathsf{F} \psi = \tv \text{ iff } \varphi = \tv \text{ and } \psi = \tv, \quad  \varphi \to_\mathsf{F} \psi = \bv \text{ iff } \varphi \neq \bv \text{ and } \psi = \bv.$$
Let us turn to the twist representation theorem. In this case, the twist-structure looks as follows \cite{wollic}:

\begin{definition}
\label{Ftwist}
    Consider a generalized Boolean algebra $\*B$. We define the \emph{full F-twist-algebra over $\*B$} as the algebra $\*B^{\bowtie}$ with universe $\{ (a_1 , a_2) \in B \times B \, | \,  a_1 \lor a_2 = 1 
\}$ and operations
\begin{itemize}
    \setlength \itemsep{0pt}   
    \item[i.] $(x_1, x_2) \land (y_1, y_2) := (x_1 \land y_1, x_2 \lor y_2)$,
    \item[ii.] $\neg(x_1, x_2) := (x_1, x_2)$, 
    \item[iii.] $\uv :=  (1, 1)$,
    \item[iv.] $(x_1, x_2) \to (y_1, y_2) := (x_1 \to y_1, x_2 \lor  y_2)$.
\end{itemize}
A \emph{F-twist-algebra over} $\*B$ is any subalgebra $\*A \leqslant \*B^{\bowtie}$ such that $\pi_1[A] = B$.
\end{definition}
A similar reasoning as in Fact \ref{notfullDF} implies that: 
\begin{fact}
\label{notfullF}
    A F-twist-algebra over $\*B$ is not necessarily full.

    \begin{proof}
        Indeed, let $\*B$ be the four-element Boolean algebra with universe $0 < a, b < 1$. Then, the set $B^{\bowtie} - \{ (a,  b) , (b, a) \}$ is closed under all the F-twist operations, thus forming the universe of a non-full F-twist-algebra, as desired.
    \end{proof}
\end{fact}
On the other hand, it is clear that every F-twist-algebra belongs to the following class: 

\begin{definition}
\label{Falg}
    An \emph{F-algebra} $(A; \land, \lor, \to, \neg, \uv)$ is an algebra of type $\la  2, 2, 2, 1, 0\ra $ such that $(A;  \land, \lor, \neg, \uv)$ is a centered Kleene algebra and the following equations are satisfied:
    \begin{enumerate}[(F1)]
\item  $x \land y = x \land (x \to y) $,
\item   $x \to y  \leq (x \land y) \lor \uv$,
\item  $(x \land y) \to z = x \to (y \to z)$,
\item  $\uv \leq x \to (y \to y)$,
\item $ \uv \leq ((x \to y) \to x) \to x $, 
\item $\Diamond (x \to y ) = \Diamond x \to \Diamond y$,
\end{enumerate}
where $\Diamond x := x \land \uv$.
\end{definition}
As in the DF case, given a F-algebra $(A; \land, \lor, \to, \neg, \uv)$, define $\Diamond(\*A):=(\Diamond[A];  \land, \lor, \to, \neg,\uv)$ by restricting the operations. We already know, by Lemma \ref{DFlemma}, that $\Diamond(\*A)$ has a reduct that is an upper-bounded distributive lattice.

\begin{lemma}
\label{Flemma}
Given a Farrell algebra $\*A$, it holds that $\Diamond(\*A)$ is a generalized Boolean algebra.
\end{lemma}
\begin{proof}
    Note that, for $a,b \in \Diamond[A]$, we have $a \leq b$ if and only if
    $a \to b = \uv$. Indeed, if $a \leq b$, then
    $a \to b = (a \land b) \to b = a \to (b \to b) \geq \uv$ by (F3).
    Conversely, if $a \to b = \uv$, then $a \land b = a \land (a \to b) = a \land \uv = a$, by (F1).
    Then, the equation $(x \land y) \to y = x \to (y \to z)$, that is, (F3), immediately entails that 
    $\land$ and $\to$ form a residual pair in $\Diamond[A]$. 
    It remains to verify Peirce's law, which is guaranteed by (F5), since $ \Diamond (((x \to y) \to x) \to x) = \uv$ holds.
\end{proof}

\begin{theorem}[F-twist representation]
\label{Ftwistrep}
     Every F-algebra $\*A$ is isomorphic to a F-twist-algebra over $\Diamond(\*A)$, as witnessed by the map $\iota \colon \*A \to \Diamond (\*A) \times \Diamond (\*A): a \mapsto (\Diamond a, \Diamond \nnot a)$.
\end{theorem}
\begin{proof}
We know from Theorem \ref{DFtwistrep}
that  $\iota$ is injective and preserves  the Kleene algebra operations. It remains to prove that it preserves the implication. Recall that $\Diamond$ preserves joins, and note that (F1) and (F2)
entail that
$(x \land y) \lor \uv = (x \to y) \lor \uv$, so
$\nnot ((x \land y) \lor \uv) =\Diamond \nnot (x \land y) = \Diamond \nnot (x \to y) = \nnot ((x \to y) \lor \uv)$, that is, 
$$\Diamond \nnot (x \land y) = \Diamond \nnot (x \to y)$$
holds. We then have:
\begin{align*}
\iota (a \to b) 
& = ( \Diamond( a \to b), \Diamond( \nnot (a \to b) ))\\
& = ( \Diamond a \to \Diamond b, \Diamond \nnot (a \land b) )& \Diamond \nnot (x \land y) = \Diamond \nnot (x \to y) \\
& =  ( \Diamond a \to \Diamond b,  \Diamond ( \nnot a  \lor \nnot b ) ) & \text{by~De Morgan's law} \\
& =  ( \Diamond a \to \Diamond b,  \Diamond \nnot a \lor \Diamond \nnot b) )& \Diamond (x \lor y ) = \Diamond x \lor \Diamond y \\
& =  ( \Diamond a, \Diamond \nnot a  )\to  (  \Diamond b,  \Diamond \nnot b  ) \\
& = \iota (a) \to \iota (b).
\end{align*}
\end{proof}

\subsection*{Cantwell}

Cantwell's intuitions in \cite{Ca08} derive from the way in which three-valued negation is introduced, namely, one sets:
$$\neg \varphi = \tv \text{ iff } \varphi = \bv, \quad \neg \varphi = \bv \text{ iff } \varphi = \tv,$$
where the remaining case is decided as $\uv$. Cantwell understands this assignment as a \emph{conditional negation}, that is, he reads $\neg \varphi$ as `$\varphi$ is false if it has a truth-value' \cite{Ca08} and then introduces an implication operator (Cooper's!) in order to study how it behaves in relation with such negation. This analysis yields the system CN presented above. 

The paper~\cite{FaOd23} studies a family of connexive systems that extend Wansing's logic C,
to which CN also belongs. These logics are all algebraizable~\cite[Thm. 50]{FaOd23},
and their algebraic models can be represented as twist-structures \cite[Def. 3]{FaOd23}.
Here, following \cite{wollic}, 
we propose an alternative presentation for the algebraic counterpart of CN, as well as an alternative twist representation. This is based on the observation that, since CN may be viewed as a  conservative expansion of De Finetti's logic by a new implication connective,
one would also expect its algebraic models to be De Finetti algebras enriched with a new operation.

\begin{definition}
\label{CNalg}
    A \emph{CN-algebra} is an algebra $\*A := (A; \land, 
    \to, \nnot, \uv)$ such that the reduct $(A; \land, 
        \nnot, \uv)$ is a centered Kleene lattice (with $\lor$ defined, as usual, through De Morgan's law) and the following equations hold:
    \begin{enumerate}[(CN1)]
        \item  $(x \land y) \to z = x \to (y \to z)$,
        \item  $\uv \leq x \to (y \to y)$,
        \item $ \uv \leq ((x \to y) \to x) \to x $,
        \item $\Diamond (x \to y ) = \Diamond x \to \Diamond y$,
        \item  $\nnot (x \to y ) =  x \to \nnot y$,
        \item 
        $\Diamond (a \land b) = \Diamond (a \land (a \to b))$,
\end{enumerate}
where $\Diamond x := x \land \uv$.
\end{definition}

Let us now look at the corresponding twist construction. 

\begin{definition}
\label{CNtwist}
    Given a generalized Boolean algebra $\*B$, 
    the \emph{full CN-twist-algebra over $\*B$} is the algebra $\*B^{\bowtie}$ with universe $B^{\bowtie} := \{ (a_1 , a_2) \in B \times B \, | \,  a_1 \lor a_2 = 1 
\}$ and operations
\begin{itemize}
    \setlength \itemsep{0pt}   
    \item[i.] $(x_1, x_2) \land (y_1, y_2) := (x_1 \land y_1, x_2 \lor y_2)$,
    \item[ii.] $\neg(x_1, x_2) := (x_2, x_1)$, 
    \item[iii.] $\uv :=  (1, 1)$,
    \item[iv.] $(x_1, x_2) \to (y_1, y_2) := (x_1 \to y_1, x_1 \to y_2)$.
\end{itemize}
We define $x \lor y : = \nnot (\nnot x \land \nnot y)$. A \emph{CN-twist-algebra over} $\*B$ is any subalgebra $\*A \leqslant {\*B}^{\bowtie}$ such that $\pi_1[A] = B$.
\end{definition}
Checking that 
every CN-twist-algebra is a CN-algebra is routine. 
A similar observation as in Fact \ref{notfullF} can be made in the present case.
\begin{remark}
\label{notfullCN}
    A CN-twist-algebra over $\*B$ is not necessarily full.
\end{remark}

The DF case already provides some information on CN. Given a CN-algebra $\*A := (A;  \land, \to, \neg, \uv)$, define $\Diamond: A \to A: x \mapsto x \land \uv$ as before. Similarly as in the DF case, by restricting the operations we can define the algebra $\Diamond(\*A) := (\Diamond[A]; \land, \lor, \to, \uv)$. We already know, by Lemma \ref{DFlemma}, that $\Diamond(\*A)$ has a reduct that is an upper-bounded distributive lattice. 
\begin{lemma}
    \label{CNlemma}
    Given a CN-algebra $\*A$, it holds that $\Diamond(\*A)$ is a generalized Boolean algebra.
\end{lemma}
\begin{proof}
   We claim that for $a,b \in \Diamond[A]$, we have $a \leq b$ if and only if
    $a \to b = \uv$. Indeed, if $a \leq b$, then $a \land b = a$, so
    $a \to b = (a \land b) \to b = a \to (b \to b) \geq \uv$ by (CN1) and (CN2). Then, since $a, b \in \Diamond (A)$, $a \to b = \Diamond a \to \Diamond b = \Diamond (a \to b) \leqslant \uv$, by (CN4) and the definition of $\Diamond[A]$.
    Conversely, if $a \to b = \uv$ then, since $a, b \in \Diamond (A)$, $a \land b = \Diamond a \land \Diamond b = \Diamond (a \land b) = \Diamond (a \land (a \to b)) = \Diamond (a \land \uv) = \Diamond \Diamond a = \Diamond a = a$ by (CN6), so $a \leqslant b$.
    Then, the equation (CN1) immediately entails that 
    $\land, \to$ form a residual pair in $\Diamond[A]$. 
    It remains to verify Peirce's law, which is guaranteed by (CN3), 
    that is, $ \Diamond (((x \to y) \to x) \to x) = \uv$.
    
\end{proof}

\begin{theorem}[CN-twist representation]
    \label{CNtwistrep}
    Every CN-algebra $\*A$ is isomorphic to a CN-twist-algebra over $\Diamond(\*A)$, as witnessed by the map $\iota \colon \*A \to \Diamond (\*A) \times \Diamond (\*A): a \mapsto (\Diamond a, \Diamond \nnot a)$.
\end{theorem}
\begin{proof}
   By Theorem \ref{DFtwistrep}, we only need to check that
that  $\iota$  preserves  the CN implication.
We then have:
\begin{align*}
\iota (a \to b) 
& = (\Diamond( a \to b), \Diamond( \nnot (a \to b) ))\\
& = (\Diamond( a \to b), \Diamond  (a \to \nnot b) )
& \text{by~(CN5)} \\
& =  (\Diamond a \to \Diamond b,  \Diamond   a  \to \Diamond \nnot b ) ) & \text{by~(CN4)} \\
& =  (\Diamond a \to \Diamond b,  \Diamond a \to \Diamond \nnot b )\\
& =  (\Diamond a, \Diamond \nnot a  )\to  ( \Diamond b,  \Diamond \nnot b  )\\
& = \iota (a) \to \iota (b).
\end{align*}
\end{proof}

The twist representation  makes it easy to compare 
the algebraic counterparts of the logics discussed so far.
As mentioned above, every CN-algebra may be viewed as a centered Kleene lattice enriched with 
a new (non-term-definable) operation $\to$ that must be given, on twist-algebras,
as per Definition~\ref{CNtwist}. In turn, 
OL-algebras may be viewed as CN-algebras enriched with a new operation 
$*$ required to satisfy the properties listed below.
Since the designated elements are in all cases $D = \{\uv, \tv\}$,
these relation are perfectly mirrored on the logical side:
by expanding Priest's LP with the constant $\uv$ we obtain De Finetti's logic,
and further expanding the latter with the CN/OL implication we obtain CN. 
As observed in~\cite{wollic}, the language of OL allows one to recover, among others,
the Kleene lattice connectives (see also~\cite[\S 4]{Riv2x} and \cite{wollic}):
\begin{align*}
    x \land_\mathsf{K} y & := (x\to_\mathsf{OL} x) \lor_\mathsf{OL} ((x \to_\mathsf{OL} y)\to_\mathsf{OL} y), \\
x \lor_\mathsf{K} y & := \neg (\neg x \land_\mathsf{K} \neg y).
\end{align*}

In order to obtain from CN the full logic OL, however, one needs to introduce at least one more
non-term-definable operation, as shown below.
\begin{remark}
\label{OLalgCNalg}
    OL-algebras can be defined as CN-algebras expanded with an operation $*$ corresponding to $\land_\mathsf{OL}$ verifying: 
\begin{enumerate}[($*$1)]
    \item $\Diamond (x * y ) = \Diamond (x \land y)$,
    \item $\Diamond \nnot (x * y ) =  (x \to \Diamond \nnot y) \land  (y \to \Diamond \nnot x)$.
\end{enumerate}

Note that, given an OL-algebra $(A; \land, \lor, \neg, \to, \uv, *)$ seen in this way, we already know by the comments above that $(\Diamond[A]; \land, \lor, \to, \uv)$ is a generalized Boolean algebra. Hence, it follows that the proof for Theorem \ref{CNtwistrep} also applies for OL; one only needs to check that $*$ is preserved by $\iota$. On one hand we have that, by using ($*$1) and ($*$2),
$$\iota(a * b) = (\Diamond (a * b), \Diamond \neg (a * b)) = (\Diamond (a \land b),  (a \to \Diamond \nnot b) \land  (b \to \Diamond \nnot a))$$
and, on the other, 
$$\iota(a)*\iota(b) = (\Diamond a \land \Diamond b, (\Diamond a \to \Diamond \neg b) \land (\Diamond b \to \Diamond \neg a)),$$
so we just need one additional axiom in order to achieve the equality: 
\begin{enumerate}[($*$3)]
    \item $\Diamond x \to \Diamond y = x \to \Diamond y$.
\end{enumerate}
\end{remark}

\section{From three to four values}

In the previous section we have 
seen how twist-algebras can be used in order to provide algebraic semantics of 
the logics of interest. 
Moreover, as the case of Cooper's logic has shown, they can also provide 
insight on the structure of the algebraic counterparts corresponding to these systems. 
In the present section we take 
the opposite direction, namely, we generalize the twist constructions and study the resulting logics. 
\\

We are going to  extend the three-valued setting
by adding a fourth truth-value
$\bot$. We can informally understand the factors of twist-algebras as follows: the first can be viewed as a \emph{truth-sensor}, and the second as a  \emph{falsity-sensor}. The overall strategy carried out in the twist semantics 
consists in identifying $\bv$ with $(0, 1)$, $\uv$ 
with $(1, 1)$ and $\tv$ with $(1, 0)$. In other words, $\bv$ is identified with `being detected by the falsity-sensor (only)', $\tv$ with `being detected by the truth-sensor (only)' and $\top$ with `being detected by both'. 
However, the twist operations 
also admit as argument the pair $(0, 0)$, which would correspond to `being detected neither by the truth nor the falsity sensors'. This pair can therefore be taken as the new truth-value $\bot$ through the previous identification: we can simply compute the tables for $(0, 0)$ in order to obtain new algebras and, as a result, new logics induced by these. In what follows, we let $A_4 := \{\bv, \bot, \top, \tv\}$.

\subsection{Adding a semantic gap}


As mentioned 
in the introduction, the truth-value $\top$ was meant to represent 
a semantic gap, that is, a situation from natural language in which a conditional assertion seems to be neither true nor false in a classical sense. 
Consider the standard examples of \emph{the Honest} and \emph{the Liar} (see, e.g., \cite{szmuc}): 
$$\text{`this statement is true'} \quad \text{ and } \quad \text{`this statement is false'.}$$
One could argue that the first statement fails to be either true or false, since we can \emph{ad libitum} assign both truth-values to it. The second, likewise, 
may be taken to be simultaneously true and false, for each time 
we try to give it a classical truth-value, it also receives the other.
Statements of the former kind are usually called \emph{hypodoxical}, those of the latter are named \emph{paradoxical} (see, e.g., \cite{priest}). 
In these terms, our formal setting has so far  only taken  under consideration hypodoxical statements;
but we could use it to treat
paradoxical statements as well. In this intuitive sense, we may speak of both semantic gaps and \textit{gluts}: if $\top$ was introduced in order to deal with `under-informative' statements, perhaps $\bot$ can model `over-informative' situations.
\\

However, the `sensor interpretation' of the four values given above, together with the fact that DF is essentially related to Priest's logic of paradox (see the preceding section), seems to indicate that $\top$ better corresponds to a formalization of semantic gluts, at least informally. Moreover, if one compares the resulting tables below with those from Belnap's logic, one clearly sees that $\top$ and $\bot$ can be identified with $\mathbf{b}$ and $\mathbf{n}$, respectively. In other words, it seems that our initial philosophical insights should be refined: $\top$ has a \textit{glutty} behavior while $\bot$ has a \textit{gappy} one. This seems to be consistent with the fact that we designate $\top$ instead of $\bot$, as it  is done in Belnap's logic. Nevertheless, these considerations stand in need of further elaboration and analysis. 
\\

Even if, at a first glance, $\top$ and $\bot$ could be understood as dealing with orthogonal phenomena from natural speech (say, two different yet parallel kinds of semantic failure), the previous observations, as well as our formal setting, will prevent them from behaving in the same way: the only exception is the case of negation, as one can check below. It is worth noting that, as consequence, one cannot hope to obtain a full three-valued fragment with $\bot$ as a logic implicit to, say, these parallel linguistic situations: the set $\{\bv, \bot, \tv\}$ fails to be closed precisely under $\to$. 

\paragraph{De Finetti.} We work out this case first because it will provide us with some insight for the CN one, as happened in the three-valued setting. The tables for the corresponding operations are the following ones: 

\begin{displaymath}
\begin{array}{|c|c|}
 & \neg_\mathsf{g} \\ 
\hline 
\bv & \tv\\
\bot & \bot\\
\top & \top \\
\tv & \bv\\
\end{array}
\quad 
\begin{array}{|c|c c c c|}
\to_\mathsf{DF}^4 & \bv & \bot & \top & \tv\\ 
\hline 
\bv & \top & \top & \top & \top \\
\bot & \bv & \bot & \bv & \bot\\
\top & \top & \top & \top & \top\\
\tv & \bv & \bot& \top & \tv \\
\end{array}
\quad 
\begin{array}{|c|c c c c|}
\land_\mathsf{K}^4 & \bv & \bot & \top & \tv\\ 
\hline 
\bv & \bv & \bv& \bv & \bv\\
\bot & \bv & \bot & \bv & \bot\\
\top & \bv & \bv & \top & \top\\
\tv & \bv & \bot& \top & \tv \\
\end{array}
\quad 
\begin{array}{|c|c c c c|}
\lor_\mathsf{K}^4 & \bv & \bot & \top & \tv\\ 
\hline 
\bv &  \bv &  \bot & \top & \tv\\
\bot & \bot & \bot & \tv & \tv \\
\top & \top & \tv &\top & \tv\\
\tv & \tv & \tv & \tv & \tv \\
\end{array}
\end{displaymath}

The implication $\to_\mathsf{DF}$ is computed by the term that defines it in the three-valued case. In this way, we may define the logic induced by $\mathrm{DFg}$ induced by the matrix $(\mathbf{DF^g_4}, \{\top, \tv\})$, where $\mathbf{DF^g_4} := (A_4; \land_\mathsf{K}^4, \lor_\mathsf{K}^4, \to_\mathsf{DF}^4, \neg_\mathsf{g})$. As we have seen before, we can simply work with $\land_\mathsf{K}^4, \lor_\mathsf{K}^4$ and $\top$. 
\\

A reasonable choice for the twist-structures in this case is the following: 
\begin{definition}
\label{DFgtwist}
     Let $\*L := (L; \land, \lor, 0, 1)$ be a bounded distributive lattice. The \emph{full DFg-twist-algebra over $\*L$} consists in an algebra ${\*L}^{\bowtie}$ with universe $L\times L$ and operations
\begin{itemize}
    \setlength \itemsep{0pt}   
    \item[i.] $(x_1, x_2) \land (y_1, y_2) := (x_1 \land y_1, x_2 \lor y_2)$,
    \item[ii.] $\neg(x_1, x_2) := (x_2, x_1)$, 
    \item[iii.] $\top :=  (1, 1)$ and $\bot := (0, 0)$.
\end{itemize}
A \emph{DFg-twist-algebra over $\*L$} is any subalgebra $\*A\leqslant {\*L}^{\bowtie}$ such that $\pi_1[A] = L$. 
\end{definition}
Obviously, the canonical map above guarantees that $\mathbf{DF_4^g}$ is isomorphic to the full DFg-twist-algebra over the two-element distributive lattice. 
\\

Now, in order to prove a twist representation result we first need an ambient class of algebras. In this case, we may simply let: 

\begin{definition}
    A \emph{DFg-algebra} is an algebra $(A; \land, \lor, \neg, 0, 1, \bot, \top)$ 
    such that the reduct $(A; \land, \lor, \neg, 0, 1 )$ is a De Morgan algebra
    and
$\bot \land \top = 0$. 
\end{definition}
Then, clearly, every DFg-twist-algebra is a DFg-algebra ($\bv$, through the canonical map, is the lower bound). Given a DFg-algebra $\*A := (A; \land, \lor, \neg, 0, \bot, \top)$, define the map $\Diamond: A \to A: x \mapsto x \land \top$. We define the algebra $\Diamond({\*A}) := (\Diamond[A]; \land, \lor, 0, \top)$, in which the operations are defined as restrictions of the ones from $\*A$. Then, the following is clear by the DF case and noting that $0$ is the desired lower bound.

\begin{lemma}
\label{DFglemma}
    Given a DFg-algebra $\*A$, it holds that $(\Diamond({\*A}); \land, \lor, 0, \top)$ is a bounded distributive lattice. 
\end{lemma}

Note that in Theorem \ref{DFtwistrep} we did not use the Kleene condition save for proving that the target universe was of the desired form. Hence, the same proof will also work for the present case.

\begin{theorem}[DFg-twist representation]
\label{DFgtwistrep}
Every DFg-algebra $\*A$ is isomorphic to a DFg-twist-algebra over $\Diamond({\*A})$, as witnessed by the map  $\iota \colon \*A \to \Diamond (\*A) \times \Diamond (\*A): a \mapsto (\Diamond a, \Diamond \nnot a)$.
\end{theorem}
\begin{proof}
    By Theorem \ref{DFtwistrep}, one only needs to prove that:
    $$\iota(\bot) = (\Diamond\bot, \Diamond\neg \bot) = (\bot \land \top, \neg \bot \land \top) =(\bot \land \top, \bot \land \top) = (0, 0).
    $$
\end{proof}

\paragraph{Cantwell.} 
We define the logic $\mathrm{CNg}$ as induced by the matrix $(\mathbf{CN^g_4}, \{\top, \tv\})$, where $\mathbf{CN^g_4} := (A_4;  \land_\mathsf{K}^4, \lor_\mathsf{K}^4, \to_\mathsf{OL}^4, \neg_\mathsf{g})$. This logic has already been considered in \cite[p. 61]{FaOd23}. Not only this but, as proven in \cite{WanConnexive05}, it turns out that the algebraizable logic MC \cite[\S 5.2]{FaOd23} coincides with $\mathrm{CNg}$. The algebraic counterpart of MC consists in class of \emph{MC-algebras} \cite[Def. 2, Def. 7]{FaOd23}.
\begin{displaymath}
\begin{array}{|c|c|}
 & \neg_\mathsf{g} \\ 
\hline 
\bv & \tv\\
\bot & \bot\\
\top & \top \\
\tv & \bv\\
\end{array}
\quad 
\begin{array}{|c|c c c c|}
\to_\mathsf{OL}^4 & \bv & \bot & \top & \tv\\ 
\hline 
\bv & \top & \top & \top & \top\\
\bot & \top & \top & \top & \top\\
\top & \bv & \bot & \top & \tv\\
\tv & \bv & \bot & \top & \tv \\
\end{array}
\quad 
\begin{array}{|c|c c c c|}
\land_\mathsf{K}^4 & \bv & \bot & \top & \tv\\ 
\hline 
\bv & \bv & \bv& \bv & \bv\\
\bot & \bv & \bot & \bv & \bot\\
\top & \bv & \bv & \top & \top\\
\tv & \bv & \bot& \top & \tv \\
\end{array}
\quad 
\begin{array}{|c|c c c c|}
\lor_\mathsf{K}^4 & \bv & \bot & \top & \tv\\ 
\hline 
\bv &  \bv &  \bot & \top & \tv\\
\bot & \bot & \bot & \tv & \tv \\
\top & \top & \tv &\top & \tv\\
\tv & \tv & \tv & \tv & \tv \\
\end{array}
\end{displaymath}

It comes in handy to recall that $\mathbf{CN^g_4}$ is in fact isomorphic to the \textit{full C-twist-algebra} \cite[Def. 3]{FaOd23} over the two-element Hilbert lattice (or, equivalently, the two-element Boolean algebra), which generates MC \cite[\S 5.2]{FaOd23}. This motivates defining the following (provisional) generalization of CN-twist-algebras:

\begin{definition}
\label{CNgtwist}
    Let $\*B$ be a Boolean algebra. We define the \emph{full CNg-twist-algebra over $\*B$} as the algebra ${\*B}^{\bowtie}$ with universe $B\times B$ and operations
\begin{itemize}
    \setlength \itemsep{0pt}   
    \item[i.] $(x_1, x_2) \land (y_1, y_2) := (x_1 \land y_1, x_2 \lor y_2)$,
    \item[ii.] $\neg(x_1, x_2) := (x_2, x_1)$, 
    \item[iii.] $(x_1, x_2) \to (y_1, y_2) := (x_1 \to y_1, x_1 \to y_2)$,
    \item[iv.] $\top :=  (1, 1)$ and $\bot := (0, 0)$.
\end{itemize}
As before, one can define $x \lor y : = \nnot (\nnot x \land \nnot y)$. A \emph{CNg-twist-algebra over} $\*B$ is any subalgebra $\*A \leqslant {\*B}^{\bowtie}$ such that $\pi_1[A] = B$. 
\end{definition}
Then, it is clear that $\mathbf{CN^g_4}$ is isomorphic to the full CNg-twist-algebra over the two-element Boolean algebra by the usual identification. Note that, since $\top$ is immediately definable by $\bot$ (letting $\top := \bot \to \bot$), we can omit $\top$ in the previous definition. Now, in order to prove a twist representation result it seems natural to employ what we already know from the DFg case. It is clear that every CNg-twist-algebra belongs to the following class:
\begin{definition}
\label{CNgalg}
    A \emph{CNg-algebra} $(A; \land, \lor, \to, \neg, 0, \bot, \top)$ is the expansion of a DFg-algebra $(A;  \land, \lor, \neg, 0, \bot, \top)$ with a binary operation $\to$ verifying:
     \begin{enumerate}[(CNg1)]
        \item  $(x \land y) \to z = x \to (y \to z)$,
        \item  $\uv \leq x \to (y \to y)$,
        \item $ \uv \leq ((x \to y) \to x) \to x $,
        \item $\Diamond (x \to y ) = \Diamond x \to \Diamond y$,
        \item  $\nnot (x \to y ) =  x \to \nnot y$,
        \item $\Diamond (a \land b) = \Diamond (a \land (a \to b))$,
\end{enumerate}
where $\Diamond x := x \land \uv$.
\end{definition}
Given a CNg-algebra $(A; \land, \lor, 0, \to, \neg, \bot, \top)$, define $\Diamond: A \to A: x \mapsto x \land \top$ as before. Similarly as in the DFg case, by restricting the operations we can define the algebra $\Diamond(\*A) := (\Diamond[A]; \land, \lor, \to, 0,  \top)$. We already know that $\Diamond(\*A)$ has a reduct that is a bounded distributive lattice. Moreover, the CN case allows us to prove, by simply noting that $0$ is a lower bound for the generalized Boolean algebra $(\Diamond[A]; \land, \lor, \to, \top)$ that:

\begin{lemma}
\label{CNglemma}
    Given a CNg-algebra $\*A$, it holds that $\Diamond(\*A)$ is a Boolean algebra.
\end{lemma}

As in the DFg case, one simply needs to check that $\bot$ is sent to $(0, 0)$ by the map below:
\begin{theorem}[CNg-twist representation]
\label{CNgtwistrep}
    Every CNg-algebra $\*A$ is isomorphic to a CNg-twist-algebra over $\Diamond(\*A)$, as witnessed by the map $\iota \colon \*A \to \Diamond (\*A) \times \Diamond (\*A): a \mapsto (\Diamond a, \Diamond \nnot a)$.
\end{theorem}

\paragraph{Farrell.}  Here, the twist-structure looks as follows: 

\begin{definition}
\label{Fgwist}
    Consider a Boolean algebra $\*B$. We define the \emph{full Fg-twist-algebra over $\*B$} as the algebra ${\*B}^{\bowtie}$ with universe $B \times B$ and operations
\begin{itemize}
    \setlength \itemsep{0pt}   
    \item[i.] $(x_1, x_2) \land (y_1, y_2) := (x_1 \land y_1, x_2 \lor y_2)$,
    \item[ii.] $\neg(x_1, x_2) := (x_2, x_1)$, 
    \item[iii.] $\top :=  (1, 1)$ and $\bot := (0, 0)$,
    \item[iv.] $(x_1, x_2) \to (y_1, y_2) := (x_1 \to y_1, x_2 \lor  y_2)$.
\end{itemize}
A \emph{Fg-twist-algebra over} $\*B$ is any subalgebra $\*A\leqslant {\*B}^{\bowtie}$ such that $\pi_1[A] = B$.
\end{definition}

This allows us to compute the corresponding table for the implication: 
\begin{displaymath}
\begin{array}{|c|c|}
 & \neg_\mathsf{g} \\ 
\hline 
\bv & \tv\\
\bot & \bot\\
\top & \top \\
\tv & \bv\\
\end{array}
\quad 
\begin{array}{|c|c c c c|}
\to_\mathsf{F}^4 & \bv & \bot & \top & \tv\\ 
\hline 
\bv & \top & \top & \top & \top\\
\bot & \top & \tv & \top & \tv\\
\top & \bv & \bv & \top & \top\\
\tv & \bv & \bot & \top & \tv \\
\end{array}
\quad\begin{array}{|c|c c c c|}
\land_\mathsf{K}^4 & \bv & \bot & \top & \tv\\ 
\hline 
\bv & \bv & \bv& \bv & \bv\\
\bot & \bv & \bot & \bv & \bot\\
\top & \bv & \bv & \top & \top\\
\tv & \bv & \bot& \top & \tv \\
\end{array}
\quad 
\begin{array}{|c|c c c c|}
\lor_\mathsf{K}^4 & \bv & \bot & \top & \tv\\ 
\hline 
\bv &  \bv &  \bot & \top & \tv\\
\bot & \bot & \bot & \tv & \tv \\
\top & \top & \tv &\top & \tv\\
\tv & \tv & \tv & \tv & \tv \\
\end{array}
\end{displaymath}
We then define the logic Fg as induced by the matrix $(\mathbf{F^g_4}, \{\top, \tv\})$, where $\mathbf{F^g_4} := (A_4; \land_\mathsf{K}^4, \lor_\mathsf{K}^4, \to_\mathsf{F}^4, \neg_\mathsf{g})$. Obviously, the preceding remark on the inter-definability of $\to_\mathsf{OL}$ and $\to_\mathsf{F}$ still holds, so we can simply refer to the CNg case. However, for the sake of completeness, let us point out that one can define: 

\begin{definition}
    A \emph{Fg-algebra} $(A; \land, \lor, \neg, \to, 0, \bot, \top)$ is the expansion of a DFg-algebra $(A;  \land, \lor, \neg, 0, \bot, \top)$ with a binary operation $\to$ verifying:
        \begin{enumerate}[(Fg1)]
\item  $x \land y = x \land (x \to y) $,
\item   $x \to y  \leq (x \land y) \lor \uv$,
\item  $(x \land y) \to z = x \to (y \to z)$,
\item  $\uv \leq x \to (y \to y)$,
\item $ \uv \leq ((x \to y) \to x) \to x $, 
\item $\Diamond (x \to y ) = \Diamond x \to \Diamond y$,
\end{enumerate}
where $\Diamond x := x \land \top$.
\end{definition}

So that every Fg-twist-algebra is in fact a Fg-algebra. Moreover, given a Fg-algebra $(A; \land, \lor, \to, \neg, 0,  \bot, \top)$, we can define $\Diamond: A \to A: x \mapsto x \land \top$ as before. Similarly as in the CNg case, by restricting the operations we can define the algebra $\Diamond(\*A) := (\Diamond[A]; \land, \lor, \to, 0,  \top)$ and, again, we already know that $\Diamond(\*A)$ has a reduct that is a bounded distributive lattice. Moreover, the F case allows us to prove, by simply noting that $0$ is a lower bound for the generalized Boolean algebra $(\Diamond[A]; \land, \lor, \to, \top)$ that:

\begin{lemma}
\label{Fglemma}
    Given a Fg-algebra $\*A$, it holds that $\Diamond(\*A)$ is a Boolean algebra.
\end{lemma}

As in the CNg case, one simply needs to check that $\bot$ is sent to $(0, 0)$ by the map below:

\begin{theorem}[Fg-twist representation]
\label{Fgtwistrep}
    Every Fg-algebra $\*A$ is isomorphic to a Fg-twist-algebra over $\Diamond(\*A)$, as witnessed by the map $\iota \colon \*A \to \Diamond (\*A) \times \Diamond (\*A): a \mapsto (\Diamond a, \Diamond \nnot a)$.
\end{theorem}

\paragraph{Cooper.} One can check that the tables for the expanded operations look as follows:

\begin{displaymath}
\begin{array}{|c|c|}
 & \neg_\mathsf{g} \\ 
\hline 
\bv & \tv\\
\bot & \bot\\
\top & \top \\
\tv & \bv\\
\end{array}
\quad 
\begin{array}{|c|c c c c|}
\to_\mathsf{OL}^4 & \bv & \bot & \top & \tv\\ 
\hline 
\bv & \top & \top & \top & \top\\
\bot & \top & \top & \top & \top\\
\top & \bv & \bot & \top & \tv\\
\tv & \bv & \bot & \top & \tv \\
\end{array}
\quad 
\begin{array}{|c|c c c c|}
\land_\mathsf{OL}^4 & \bv & \bot & \top & \tv\\ 
\hline 
\bv & \bv & \bv& \bv & \bv\\
\bot & \bv & \bv & \bot & \bot\\
\top & \bv & \bot & \top & \tv\\
\tv & \bv & \bot& \tv & \tv \\
\end{array}
\quad 
\begin{array}{|c|c c c c|}
\lor_\mathsf{OL}^4 & \bv & \bot & \top & \tv\\ 
\hline 
\bv &  \bv &  \bot & \bv & \tv\\
\bot & \bot & \tv & \bot & \tv \\
\top & \bv & \bot &\top & \tv\\
\tv & \tv & \tv & \tv & \tv \\
\end{array}
\end{displaymath}
It is worth noting that the implication operation coincides with the one of \textit{material connexive logic} \cite[Rmk. 4.3]{OmWa20}. We can then consider the logic $\mathrm{OLg}$ induced by the matrix $(\mathbf{OL^g_4}, \{\top, \tv\})$, where $\mathbf{OL^g_4} := (A_4; \land_\mathsf{OL}^4, \lor_\mathsf{OL}^4, \to_\mathsf{OL}^4,\neg_\mathsf{g})$. Here one would naturally think of algebraic models as OL-twist-algebras with an additional fixpoint for negation. The problem is that one needs $(0, 0)$ to belong to the corresponding algebra, so this implies relaxing the former conditions imposed on the universe. One could simply consider:
\begin{definition}
\label{OLgtwist}
     Let $\*B$ be a Boolean algebra. The \emph{full OLg-twist-algebra over $\*B$} is the algebra ${\*B}^{\bowtie}$ consisting of universe $B\times B$ and operations
\begin{itemize}
    \setlength \itemsep{0pt}   
    \item[i.] $(x_1, x_2) * (y_1, y_2) := (x_1 \land y_1, (x_1 \to y_2)\land (y_1 \to x_2)),$ 
    \item[ii.] $(x_1,x_2) \to (y_1, y_2) := (x_1 \to y_1, x_1 \to y_2),$ 
    \item[iii.] $\neg (x_1, x_2) := (x_2, x_1)$.
\end{itemize}
An \emph{OLg-twist-algebra over} $\*B$ is any subalgebra $\*A\leqslant {\*B}^{\bowtie}$ such that $\pi_1[A] = B$.
\end{definition}


Then, one can check that $\mathbf{OL^g_4}$ is isomorphic to the OLg-twist-algebra over the two-element (generalized) Boolean algebra $\{0, 1\}$ via the identification mentioned at the beginning of this section. Hence, the task here is to replicate the twist representation result that we have presented above \cite{Riv2x}.
Now, since OL-algebras are CN-algebras expanded with an additional operation corresponding to $\land_\mathsf{OL}$ (see Remark \ref{OLalgCNalg} above), it seems natural to expand Definitions \ref{CNgtwist} and \ref{CNgalg} accordingly\footnote{In general, one can replicate the construction given in \cite{FaOd23} for N4-lattices, where the corresponding twist-structure is not necessarily full but the elements are pairs $(a, b)$ verifying that $a \land b \in \Delta$ and $a \lor b \in \nabla$, for some ideal $\Delta \subseteq B$ and some filter $\nabla \subseteq B$ (see \cite{FaOd23} for further details).}.

\begin{definition}
\label{OLgtwist2}
    Let $\*B := (B; \land, \lor, \to, \neg, 0, 1)$ be a Boolean algebra. We define the \emph{full OLg-twist-algebra over $\*B$} as the full CNg-twist-algebra ${\*B}^{\bowtie}$ over $\*B$ expanded with the operation:
$$(x_1, x_2) * (y_1, y_2) := (x_1 \land y_1, (x_1 \to y_2)\land (y_1 \to x_2)).$$
An \emph{OLg-twist-algebra over} $\*B$ is any subalgebra $\*A \leqslant {\*B}^{\bowtie}$ such that $\pi_1[A] = B$. 
\end{definition}

\begin{definition}
    An \emph{OLg-algebra} $(A; \land, \lor, \to, *, \neg, 0,  \bot, \top)$ is a CNg-algebra $(A; \land, \lor, \to,\neg,  0, \bot, \top)$ expanded with an operation $*$ verifying ($*$1)-($*$3) from Remark \ref{OLalgCNalg}, namely,
    \begin{enumerate}[($*$1)]
    \item $\Diamond (x * y ) = \Diamond (x \land y)$,
    \item $\Diamond \nnot (x * y ) =  (x \to \Diamond \nnot y) \land  (y \to \Diamond \nnot x)$,
    \item $\Diamond x \to \Diamond y = x \to \Diamond y$.
\end{enumerate}
\end{definition}

Given an OLg-algebra $(A; \land, \lor, \to, *, \neg, 0, \bot, \top)$, let us set $\Diamond: A \to A: x \mapsto \Diamond x := x \land \top$ as in the CNg case. Similarly, by restricting the operations we can define the algebra $(\Diamond[A];  \land, \lor, \to, 0, \top)$, which has already the structure of Boolean algebra by Lemma \ref{CNglemma}. Thus, we only need to check that the map from Theorem \ref{CNgtwistrep} preserves $*$, but this was already seen in Remark \ref{OLalgCNalg}.

\begin{theorem}[OLg-twist representation]
\label{OLgtwistrep}
    Every OLg-algebra $\*A$ is isomorphic to a OLg-twist-algebra over $\Diamond(\*A)$, as witnessed by the map $\iota \colon \*A \to \Diamond (\*A) \times \Diamond (\*A): a \mapsto (\Diamond a, \Diamond \nnot a)$.
\end{theorem}

\subsection{Adding a falsity}

\paragraph{Motivation.} In the preceding cases, the truth-values $\bot$ and $\top$ behave as fixpoints for negation, which is consistent with the intuition that both differently model semantic anomalies from natural speech. There is yet another idea that can be traced back to Dummett, namely, that the new truth-values actually codify \textit{variations} of truth and falsity. In this way, $\top$ could be regarded as some kind of truth, while the new value $\bot$ would correspond to some kind of falsity. As Dummett puts it:
\begin{quote}
    I once imagined a case in which a language contained a negation operator
‘$-$’ which functioned much like our negation save that it made
‘$-(A \to B)$’ equivalent to ‘$A \to -B$’, where ‘$\to$’ is the ordinary two-valued
implication. In this case, the truth or falsity of ‘$-(A \to B)$’ would
not depend solely on the truth or falsity of ‘$A \to B$’, but on the particular
way in which ‘$A \to B$’ was true (whether by the truth of both constituents
or by the falsity of the antecedent). This would involve the use of three-valued
truth-tables, distinguishing two kinds of truth. In the same way,
it might be necessary to distinguish two kinds of falsity. \cite[p. 209]{humberstone}
\end{quote}

It is worth noting that the intuition that Dummett is vindicating is captured by the Boethius theses (see above). Even if indicative conditionals fail to be informative when their premise is false, Dummett prefers not to simply consider them as lacking a definite truth-value but, instead, as being true in some sense. Perhaps we could say, after all, that a vacuously true indicative conditional is indeed true \textit{in some sense}. In the same way, the negation of such statement would be false \textit{in some sense} or, accordingly, vacuously false. 
\\

Therefore, this approach relies on expanding the negation-free fragments of the preceding four-valued logics with a new negation. It is important to note that here we need an underlying negation operation $\neg$ in the corresponding twist-structures factors, since we wish to define $-(x, y) := (\neg x, \neg y)$. This new connective could be seen as the extension of classical negation for the non classical values\footnote{\label{footnote1}In \cite[\S 5.2]{avron} the \textit{conflation} operator is defined to fix the classical truth-values and interchange $\bot$ and $\top$. It is clear, however, that even if we could consider this operation instead, the resulting systems would not satisfy the basic philosophical requirements for a reasonable logic of indicative conditionals.}. Note that, as in the previous section, a full three-valued fragment with $\bot$ is not possible to obtain: in this case, both the negation and the implication operations prevent us from having such situation. 

\paragraph{Farrell.} Let us introduce the corresponding tables: 

\begin{displaymath}
\begin{array}{|c|c|}
 & \neg_\mathsf{f}\\ 
\hline 
\bv & \tv\\
\bot & \top\\
\top & \bot \\
\tv & \bv\\
\end{array}
\quad 
\begin{array}{|c|c c c c|}
\to_\mathsf{F}^4 & \bv & \bot & \top & \tv\\ 
\hline 
\bv & \top & \top & \top & \top\\
\bot & \top & \tv & \top & \tv\\
\top & \bv & \bv & \top & \top\\
\tv & \bv & \bot & \top & \tv \\
\end{array}
\quad 
\begin{array}{|c|c c c c|}
\land_\mathsf{K}^4 & \bv & \bot & \top & \tv\\ 
\hline 
\bv & \bv & \bv& \bv & \bv\\
\bot & \bv & \bot & \bv & \bot\\
\top & \bv & \bv & \top & \top\\
\tv & \bv & \bot& \top & \tv \\
\end{array}
\quad 
\begin{array}{|c|c c c c|}
\lor_\mathsf{K}^4 & \bv & \bot & \top & \tv\\ 
\hline 
\bv &  \bv &  \bot & \top & \tv\\
\bot & \bot & \bot & \tv & \tv \\
\top & \top & \tv &\top & \tv\\
\tv & \tv & \tv & \tv & \tv \\
\end{array}
\end{displaymath}

Hence, we define the logic Ff as induced by the matrix $(\mathbf{F_4^f}, \{\top, \tv\})$, where $\mathbf{F_4^f} := (A_4; \land_\mathsf{K}^4, \lor_\mathsf{K}^4, \to_\mathsf{F}^4, \neg_\mathsf{f})$. The preceding tables have been computed according to the next definition:

\begin{definition}
\label{Fftwist}
    Let $\*B$ and $\*D$ be a Boolean algebra and a De Morgan algebra, respectively. We define the \emph{full Ff-twist-algebra over $\*B$ and $\*D$} as the algebra $(\*B \times \*D)^{\bowtie}$ with universe $B \times D$ and operations
\begin{itemize}
    \setlength \itemsep{0pt}   
    \item[i.] $(x_1, x_2) \land (y_1, y_2) := (x_1 \land_{\*B} y_1, x_2 \lor_{\*D} y_2)$,
    \item[ii.] $\neg(x_1, x_2) := (\neg_{\*B} x_1, \neg_{\*D} x_2)$, 
    \item[iii.] $(x_1, x_2) \to (y_1, y_2) := (x_1 \to_{\*B} y_1, x_2 \lor_{\*D} y_2)$, 
    \item[iv.] $\bot := (0_{\*B}, 0_{\*D}).$
\end{itemize}
    Additionally, one can set $x \lor y := \neg ( \neg x \land \neg y)$. An \emph{Ff-twist-algebra over $\*B$ and $\*D$} is any subalgebra $\*A \leqslant(\*B \times \*D)^{\bowtie}$ such that $\pi_1[A] = B$. 
\end{definition}

In the preceding cases of CNg and Fg, starting with DFg was quite instructive, since it allowed use to first study the implication-free fragment of our full language. Here, it seems reasonable to proceed similarly. However, the change of the negation operator compromises the definability results above and, in particular, the definability of $\to_\mathsf{DF}$ in terms of the remaining operations in the language. Philosophically, it is then clear that De Finetti's logic is \emph{not} captured by having the strong Kleene operations and adding \emph{some} negation: this seems more like an accidental feature that holds in the DF and DFg cases. Moreover, and surprisingly enough, it turns out that: 

\begin{fact}
\label{fact:Ffdef}
      The term $(\top \land \neg x) \lor (x \land y)$ defines $\to$ in Ff-twist-algebras.
      
    \begin{proof}
        Indeed, one can compute: 
        $$((1, 1) \land (\neg x_1, \neg x_2)) \lor (x_1 \land y_1, x_2 \lor y_2) = (\neg x_1, 1) \lor (x_1 \land y_1, x_2 \lor y_2)$$
        and this is equal to $(\neg x_1 \lor (x_1 \land y_1), x_2 \lor y_2) = (\neg x_1 \lor y_1, x_2 \lor y_2)$. Now, the first component coincides with $x_1 \to y_1$ in Boolean algebras, as desired. 
    \end{proof}
\end{fact}

This Remark strongly suggests that we can simply drop out the implication operation and consider the constant $\top$ as primitive instead. Hence, we may define RFf (for `restricted' Ff) to be the logic induced by $(\mathbf{RF^f_4}, \{\top, \tv\})$, where $\mathbf{RF^f_4} := (A_4; \land_\mathsf{K}, \lor_\mathsf{K}, \neg_\mathsf{f}, \top)$, and then the twist construction will look as follows: 

\begin{definition}
\label{RFftwist}
     Let $\*D_1, \*D_2$ be two De Morgan algebras. The \emph{full RFf-twist-algebra over $\*D_1$ and $\*D_2$} consists in an algebra $\*D_1 \bowtie \*D_2$ with universe $\*D_1\times \*D_2$ and operations
\begin{itemize}
    \setlength \itemsep{0pt}   
    \item[i.] $(x_1, x_2) \land (y_1, y_2) := (x_1 \land_{\*D_1} y_1, x_2 \lor_{\*D_2} y_2)$,
    \item[ii.] $\neg (x_1, x_2) := (\neg_{\*D_1}\, x_1, \neg_{\*D_2}\, x_2)$, 
    \item[iii.] $\bot := (0_{\*D_1}, 0_{\*D_2})$.
\end{itemize}
Note that we may simply let $\top := \neg\bot$ and, similarly, we can define $x \lor y : = \neg (\neg x \land \neg y)$. An \emph{RFf-twist-algebra over $\*D_1$ and $\*D_2$} is any subalgebra $\*A \leqslant \*D_1 \bowtie \*D_2$ such that $\pi_1[A] = D_1$. 
\end{definition}

{
\begin{remark}
    \label{rmk:90deg}
    Having all four constants available, from the `truth order' $\land, \lor$ one can explicitly define, via the so-called \textit{90-degree lemma} \cite[Lemma 1.5]{jungriv}, the `knowledge order', usually considered on bi-lattices (see, e.g., \cite{jungriv}). The definitions are: $
x \sqcap y := (x \land \bot) \lor (y \land \bot) \lor (x \land y)
$ and $
x \sqcup y := (x \land \top) \lor (y \land \top) \lor (x \land y)
$. Additionally, let us note that in the preceding definition we could have defined $(x_1, x_2) \land (y_1, y_2) := (x_1 \land_{\*D_1} y_1, x_2 \land_{{\*D}'_2} y_2)$, where ${\*D}'_2$ is the dual of $\*D_2$. In this case, then, $\bot := (0_{\*D_1}, 1_{\*D'_2})$.
\end{remark}
}

We must note, however, that Fact \ref{fact:Ffdef} only works in case the ambient structure is a Ff-twist-algebra: in general, it is clear that Ff-twist and RFf-twist-algebras do not coincide, since the first factor in the former case is a Boolean algebra while in the latter is a De Morgan algebra. Let us prove then two different results for these cases. We will need, of course, some extra requirements for the Ff case to work properly (precisely in the point of recovering the Boolean algebra structure from the corresponding factor). Because of this, we will work out a proof for the twist representation of Ff as a slight variation of that of RFf. We begin by defining: 
\begin{definition}
    A \emph{RFf-algebra} $(D; \land, \lor, \neg, 0, 1, \top)$ is a bounded De Morgan lattice $(D; \land, \lor, \neg, 0, 1)$ such that: 
    \begin{enumerate}[(RFf1)]
        \item $\top \approx 1 \Rightarrow x \approx y$,

        \item $\top \land \neg \top = 0$, 
        \item $a \land \top = b \land \top$ and $a \lor \neg\top = b \lor \neg \top$ imply that $a = b$. 
    \end{enumerate}
    Observe that we also have $\neg 0 = 1$. We let $\bot := \neg \top$.
\end{definition}

\begin{remark}
\label{remark:quasieq}
    Note that, what (RFf1) tells us is that $\top$ (and $\neg \top$) receive different values to those of $0$ and $1$. We will also make use of this condition in what follows. 
\end{remark}

Now, given an RFf-algebra $\*D := (D; \land, \lor, \neg, 0, 1, \bot, \top)$, we set $\Diamond x := x \land \top$ and $\Box x := x \land \bot$. By restricting the operations from $\*D$, one can define
$$\Diamond(\*D) := (\Diamond[D]; \land, \lor, \neg_\Diamond, 0, \top) \text{ and } \Box(\* D) := (\Box [D]; \land, \lor, \neg_\Box, 0, \bot),$$
where $\neg_\Diamond \Diamond a := \Diamond \neg a$ and $\neg_\Box \Box a := \Box \neg a$. Then, it is clear that these algebras are De Morgan algebras: for instance, one can easily check that both $\neg_\Diamond$ and $\neg_\Box$ preserve tops and bottoms. Then, we can prove: 
\begin{theorem}[RFf-twist representation]
\label{RFftwistrep}
    Let $\*D$ be an RFf-algebra. Then, the map $\iota: \* D \to  \Diamond (\*D)\times \Box(\*D): a \mapsto (\Diamond a, \Box \neg a)$ witnesses an isomorphism between $\*D$ and the corresponding RFf-twist-structure.
\end{theorem}
\begin{proof}
We need to check the following requirements:
    \begin{itemize}
     \item[--] $\iota$ is injective. Indeed, $\iota(a) = \iota(b)$ implies $\Diamond a = \Diamond b$ and $\Box \neg a = \Box \neg b$. This means that $a \land \top = b \land \top$ and $\neg a \land \bot = \neg b \land \bot$. The latter gives us that $\neg (a \lor \top) = \neg (a \lor \top)$, so $a \lor \top = b \lor \top$, since $\neg \neg x = x$ holds. Then, by definition, we have already that $a = b$\footnote{Of course, a more elegant (logical-looking) condition to appear in the definition of RFf-algebras would be desirable, but the one stated above works provisionally. One may also wonder if this quasi-equation can be replaced by a set of reasonable equations, as in the three-valued case.}. 
    \item[--] $\iota$ preserves $\land$. On the one hand, $\iota(a \land b) = (\Diamond(a \land b), \Box \neg (a \land b))$ and, on the other, $\iota(a) \land \iota(b) = (\Diamond a \land \Diamond b, \Box \neg a \lor \Box \neg b)$. But $\Box \neg  a \lor \Box \neg b = \Box \neg (a \land b)$ holds, as desired. 
    \item[--] $\iota$ preserves $\neg$. Indeed, $\iota(\neg a) = (\Diamond \neg a, \Box \neg \neg a) =  (\Diamond \neg a, \Box a)$, where $\neg \iota(a) = (\neg_\Diamond \Diamond a, \neg_\Box \Box \neg a)$ $= (\Diamond \neg a, \Box \neg \neg a) =$$ (\Diamond \neg a, \Box a)$, as desired.
    \item[--] $\iota$ preserves $\bot$, since  $\iota(\bot) = (\Diamond \bot, \Box \neg \bot) = (\bot \land \top, \top \land \bot) = (0, 0)$.
    \end{itemize}
\end{proof}

Accordingly, we have: 

\begin{definition}
    A \emph{Ff-algebra} $(D; \land, \lor,\to,  \neg,  0, 1, \bot, \top)$ is a RFf-algebra $(D; \land, \lor, \neg, 0, 1, \bot, \top)$ expanded with an operation $\to$ such that the following conditions hold: 
    \begin{enumerate}[(Ff1)]
        \item $x \land y = x \land (x \to y)$,
        \item $(x \land y) \to z = x \to (y \to z)$,
        \item $\top \leq ((x \to y) \to x) \to x$,
        \item $\Diamond (x \to y ) = \Diamond x \to \Diamond y$, 
        \item $\Box\neg (x \to y) = \Box \neg (x \land y)$, 
    \end{enumerate}
    where $\Diamond x := x \land \top$ and $\Box x := x \land \bot$. 
\end{definition}
One can check that every Ff-twist-algebra is in fact a Ff-algebra. Now, define similarly as in the RFf case: 
$$(\Diamond[D]; \land, \lor, \to,\neg_\Diamond, 0, \top) \text{ and } (\Box[D]; \land, \lor, \neg_\Box, 0, \bot),$$
where, remember, $\neg_\Diamond \Diamond a := \Diamond \neg a$ and $\neg_\Box \Box a := \Box \neg a$. In order to prove twist representation for these algebras, we first need a preliminary observation. 
\begin{lemma}
\label{Fflemma}
    Given a Ff-algebra $\*D$, it holds that $\Diamond(\*D)$ is a Boolean algebra and that $\Box(\*D)$ is a De Morgan algebra. 
\end{lemma}
\begin{proof}
    Indeed, by the RFf case we already know that $\Box(\*D)$ is a De Morgan algebra. For $\Diamond(\*D)$ we can recall Lemma \ref{Flemma} and note that the required properties hold by definition plus the fact that $\Diamond[D]$ is bounded. 
\end{proof}

\begin{theorem}[Ff-twist representation]
\label{Fftwistrep}
    Let $\*D$ be a Ff-algebra. Then, the map $\iota: \*D \to  \Diamond (\*D)\times \Box(\*D): a \mapsto (\Diamond a, \Box \neg a)$ witnesses an isomorphism between $\*D$ and the corresponding Ff-twist-structure.
\end{theorem}
\begin{proof}
    By Theorem \ref{RFftwist} we only need to check that $\iota$ preserves $\to$. On one hand, 
    $$\iota(a \to b) = (\Diamond(a \to b), \Box \neg (a \to b)) = (\Diamond a\to \Diamond b, \Box \neg (a \land  b)),$$
    and, on the other,
    $$\iota(a) \to \iota(b) = (\Diamond a \to \Diamond b, \Box \neg a \lor \Box \neg b) = (\Diamond a \to \Diamond b, \Box \neg (a \land b)),$$
    as desired.
\end{proof}

\paragraph{De Finetti.} We define DFf as the logic induced by $(\mathbf{DF^f_4}, \{\top, \tv\})$, where $\mathbf{DF^f_4} := (A_4; \land_\mathsf{K}^4, \lor_\mathsf{K}^4, \to_\mathsf{DF}^4, \neg_\mathsf{f})$.

\begin{displaymath}
\begin{array}{|c|c|}
 & \neg_\mathsf{f} \\ 
\hline 
\bv & \tv\\
\bot & \top\\
\top & \bot \\
\tv & \bv\\
\end{array}
\quad 
\begin{array}{|c|c c c c|}
\to_\mathsf{DF}^4 & \bv & \bot & \top & \tv\\ 
\hline 
\bv & \top & \top & \top & \top \\
\bot & \bv & \bot & \bv & \bot\\
\top & \top & \top & \top & \top\\
\tv & \bv & \bot& \top & \tv \\
\end{array}
\quad 
\begin{array}{|c|c c c c|}
\land_\mathsf{K}^4 & \bv & \bot & \top & \tv\\ 
\hline 
\bv & \bv & \bv& \bv & \bv\\
\bot & \bv & \bot & \bv & \bot\\
\top & \bv & \bv & \top & \top\\
\tv & \bv & \bot& \top & \tv \\
\end{array}
\quad 
\begin{array}{|c|c c c c|}
\lor_\mathsf{K}^4 & \bv & \bot & \top & \tv\\ 
\hline 
\bv &  \bv &  \bot & \top & \tv\\
\bot & \bot & \bot & \tv & \tv \\
\top & \top & \tv &\top & \tv\\
\tv & \tv & \tv & \tv & \tv \\
\end{array}
\end{displaymath}

The underlying structure that allows the computation of the previous tables looks as follows:

\begin{definition}
\label{DFftwist}
     Let $\*{D}_1, \*{D}_2$ be two De Morgan algebras. Let $\rho: \*{D}_2 \xhookrightarrow{} \*{D}_1$ be an embedding. The \emph{full DFf-twist-algebra over $\*{D}_1$ and $\*{D}_2$} consists in an algebra $\*{D}_1 \bowtie \*{D}_2$ with universe $D_1\times D_2$ and operations
\begin{itemize}
    \setlength \itemsep{0pt}   
    \item[i.] $(x_1, x_2) \land (y_1, y_2) := (x_1 \land_{\*D_1} y_1, x_2 \lor_{\*D_2} y_2)$,
    \item[ii.] $\neg (x_1, x_2) := (\neg_{\*D_1} x_1, \neg_{\*D_2} x_2)$, 
    \item[iii.] $(x_1, x_2) \to (y_1, y_2):= (\rho(x_2) \lor_{\*D_1} (x_1 \land_{\*D_1} y_1), x_2 \lor_{\*D_2} y_2),$
    \item[iv.] $\bot := (0_{\*D_1}, 0_{\*D_2})$.
\end{itemize}
Note that we may simply let $\top := \neg \bot$ and, similarly, we can define $x \lor y : = \neg ( \neg x \land \neg y)$. A \emph{DFf-twist-algebra over $\*{D}_1$ and $\*{D}_2$} is any subalgebra $\*A \leqslant \*{D}_1 \bowtie \*{D}_2$ such that $\pi_1[A] = D_1$. 
\end{definition}

Let us observe that, in the previous definition, proceeding as in Remark \ref{toDF} allows us to see that
$$(x_1, x_2) \to (y_1, y_2):= (\rho(x_2) \lor (x_1 \land y_1), x_2 \lor y_2) =  ((x_1, x_2) \land (y_1, y_2)) \lor (\rho(y_1), 1),$$
so, in item (iii) above, we can consider the operation $!(x, y) := (\rho(y), 1)$ instead. Thus, we can refine our previous twist construction: 

\begin{definition}
\label{DFftwist2}
    Let $\*{D}_1, \*{D}_2$ be two De Morgan algebras. Let $\rho: \*{D}_2 \xhookrightarrow{} \*{D}_1$ be an embedding. The \emph{full DFf-twist-algebra over $\*{D}_1$ and $\*{D}_2$} consists in an algebra $\*{D}_1 \bowtie \*{D}_2$ with universe $D_1\times D_2$ and operations
\begin{itemize}
    \setlength \itemsep{0pt}   
    \item[i.] $(x_1, x_2) \land (y_1, y_2) := (x_1 \land_{\*D_1} y_1, x_2 \lor_{\*D_2} y_2)$,
    \item[ii.] $\neg (x_1, x_2) := (\neg_{\*D_1} x_1, \neg_{\*D_2} x_2)$, 
    \item[iii.] $!(x_1, x_2) := (\rho(x_2), 1_{\*D_2}),$
    \item[iv.] $\bot := (0_{\*D_1}, 0_{\*D_2})$.
\end{itemize}
As before, one lets $\top := \neg \bot$, $x \lor y : = \neg ( \neg x \land \neg y)$ and $x \to y := (x \land y) \lor  \, !y$. A \emph{DFf-twist-algebra over $\*{D}_1$ and $\*{D}_2$} is any subalgebra $\*A \leqslant \*{D}_1 \bowtie \*{D}_2$ such that $\pi_1[A] = D_1$. 
\end{definition}

Let us turn now to the twist representation result. First, let us observe that every DFf-twist-algebra belongs to the following class: 

\begin{definition}
    A \emph{DFf-algebra} $(D; \land, \lor, \neg, !, 0, 1, \top)$ is a bounded De Morgan lattice $(D; \land, \lor, \neg, 0, 1)$ expanded with a unary operation $!$ verifying: 
    \begin{enumerate}[(DFf1)]
       \item $\top \approx 1 \Rightarrow x \approx y$ (see Remark \ref{remark:quasieq}),
        \item $\top \land \neg \top= 0$, 
        \item $\Box a = \Box b$ implies that $\Diamond ! a = \Diamond ! b$,
        \item $\Diamond a = \Diamond b$ and $\Box a = \Box b$ imply that $a = b$,
        \item $\Diamond ! a = \Diamond ! b$ implies $\Box a = \Box b$,
        \item $\Diamond ! (a \lor b) = \Diamond ! a \land \Diamond ! b$,
        \item $\Diamond ! \neg a = \Diamond \neg ! a$, 
        \item $!\bot = 0$ and $!0 = \top$, 
        \item $\Box ! a = 1$,
    \end{enumerate}
     where $\bot := \neg \top$, $\Diamond x :=x \lor \bot$ and $\Box x := x \land \bot$.
\end{definition}

Now, given a DFf-algebra $\*D$, define $\Diamond x :=x \lor \bot$ and $\Box x := x \land \bot$ as above. Let us set $\Diamond(\*D):=(\Diamond[D]; \land, \lor, \neg_\Diamond, \bot, 1),$
where $\neg_\Diamond \Diamond a := \Diamond \neg a$. Moreover, consider: 
$\Box(\*D) := (\Box(D); \land_\Box, \lor_\Box, \neg_\Box, 0_\Box, 1_\Box),$
where: $\Box a \land_\Box \Box b := \Box (a \lor b)$; $\Box a \lor_\Box \Box b := \Box (a \land b)$; $\neg_\Box \Box a := \Box \neg a$; $\Box a \to_\Box \Box b := \neg_\Box \Box a \lor_\Box \Box b = \Box (\neg a \land b)$; $0_\Box := \bot$ and $1_\Box := 0$. 

\begin{lemma}
    Given a DFf-algebra $\*D$, it holds that both $\Diamond(\*D)$ and $\Box(\*D)$ are De Morgan algebras. 
\end{lemma}
\begin{proof}
    It is clear that both of them are De Morgan lattices. Additionally, they are bounded by definition, so the result follows. 
\end{proof}

Note that the definition of DFf-twist-algebras forces us to explicitly provide an embedding of De Morgan algebras between the factors of our desired twist-structure. One can check that the definition of DFf-algebras allows the following map to be well defined, injective and operations-preserving:
\begin{lemma}
    Given a DFf-algebra $\*D$, the map $\rho: \Box(\*D) \to \Diamond(\*D): \Box a \mapsto \Diamond !a$ is an embedding of De Morgan algebras.
\end{lemma}
\begin{proof}
    We have to check the following conditions and verify them by means of the definitions above: 
\begin{itemize}
    \item[--] $\rho$ is well defined. If $\Box a = \Box b$ then $\rho(\Box a) = \Diamond ! a = \Diamond ! b = \rho(b)$. 
    \item[--] $\rho$ is injective. We need that, if $\rho(\Box a) = \Diamond ! a = \Diamond ! b = \rho(\Box b)$ then $\Box a = \Box b$, which holds indeed by definition.
    \item[--] $\rho$ preserves $\land$. Indeed, $\rho(\Box a \land_\Box \Box b) = \rho(\Box (a \lor b)) = \Diamond ! (a \lor b) = \Diamond ! a \land \Diamond ! b = \rho(\Box a) \land \rho(\Box b)$.
    \item[--] $\rho$ preserves $\neg$. Indeed, $\rho(\neg_\Box \Box a) = \rho(\Box \neg a) = \Diamond ! \neg a = \Diamond \neg ! a = \neg_\Diamond \Diamond ! a = \neg_\Diamond \rho(\Box a)$. 
    \item[--] $\rho$ preserves extrema. Indeed, $\rho(0_\Box) = \rho(\bot) = \rho(\Box \bot) = \Diamond ! \bot = \Diamond 0 = \bot$. On the other hand, $\rho(1_\Box) = \rho(0) = \rho(\Box 0) = \Diamond ! 0 = \Diamond \top = 1$, as desired.
\end{itemize}
\end{proof}

Now we can finally prove: 

\begin{theorem}[DFf-twist representation]
     Let $\*D$ be a DFf-algebra. Then, the map $\iota: \*D \to  \Diamond (\*D)\times \Box(\*D): a \mapsto (\Diamond a, \Box a)$ witnesses an isomorphism between $\*D$ and the corresponding DFf-twist-structure.
\end{theorem}
\begin{proof}
The proof consists, once again, in checking the following requirements one by one: 
 \begin{itemize}
     \item[--] $\iota$ is injective. This means that $\Diamond a = \Diamond b$ and $\Box a = \Box b$ imply $a = b$, which holds by the definition of DFf-algebras.
    \item[--] $\iota$ preserves $\land$. We have: $\iota(a \land b) = (\Diamond(a \land b), \Box (a \land b)) = (\Diamond a \land \Diamond b, \Box a \lor_\Box \Box b) = \iota(a) \land \iota(b),$ as we wanted.
    \item[--] $\iota$ preserves $\neg$. Indeed, $\iota(\neg a) = (\Diamond \neg a, \Box \neg a) = (\neg_\Diamond \Diamond a, \neg_\Box \Box a) = \neg \iota(a)$.
    \item[--] $\iota$ preserves $\bot$. This holds because $\iota(\bot) = (\Diamond \bot, \Box \bot) = (\bot, \bot)$, as desired.
    \item[--] $\iota$ preserves $!$. Here we have, on the one hand $\iota(!a) := (\Diamond !a, \Box !a)$ and, on the other, $!\iota(a) = !(\Diamond a, \Box a) = (\rho(\Box a), 1).$ But these are equal in virtue of the definitions above.
    \end{itemize}
\end{proof}

\paragraph{Cantwell.} We define the logic CNf as induced by $(\mathbf{CNf}_4, \{\top, \tv\})$, where one sets $\mathbf{CNf}_4 := (A_4;  \land_\mathsf{K}^4, \lor_\mathsf{K}^4, \to_\mathsf{OL}^4, \neg_\mathsf{f})$\footnote{Dummett comments above suggested unifying both a classical negation for two kinds of truth while at the same time preserving Boethius theses, and these requirements cannot be met in the case of CNf: $\neg_\mathsf{f}(x \to_\mathsf{OL} \bv)$ lies always in $\{\tv, \bot\}$ but $x \to_\mathsf{OL} \tv$ can receive value $\top$.}. 
\begin{displaymath}
\begin{array}{|c|c|}
 & \neg_\mathsf{f} \\ 
\hline 
\bv & \tv\\
\bot & \top\\
\top & \bot \\
\tv & \bv\\
\end{array}
\quad 
\begin{array}{|c|c c c c|}
\to_\mathsf{OL}^4 & \bv & \bot & \top & \tv\\ 
\hline 
\bv & \top & \top & \top & \top\\
\bot & \top & \top & \top & \top\\
\top & \bv & \bot & \top & \tv\\
\tv & \bv & \bot & \top & \tv \\
\end{array}
\quad 
\begin{array}{|c|c c c c|}
\land_\mathsf{K}^4 & \bv & \bot & \top & \tv\\ 
\hline 
\bv & \bv & \bv& \bv & \bv\\
\bot & \bv & \bot & \bv & \bot\\
\top & \bv & \bv & \top & \top\\
\tv & \bv & \bot& \top & \tv \\
\end{array}
\quad 
\begin{array}{|c|c c c c|}
\lor_\mathsf{K}^4 & \bv & \bot & \top & \tv\\ 
\hline 
\bv &  \bv &  \bot & \top & \tv\\
\bot & \bot & \bot & \tv & \tv \\
\top & \top & \tv &\top & \tv\\
\tv & \tv & \tv & \tv & \tv \\
\end{array}
\end{displaymath}

As before, it is instructive to ask on the form of CNf-twist-algebras. The preceding cases have shown how we need to allow different factor algebras in our twist representation. Not only this but, comparing with Definition \ref{DFftwist}, we also need here some map for the definition of the implication operation, since it is the only one not determined component-wise by the inner operations of the factor algebras:

\begin{definition}
    \label{CNftwist}
    Let $\*{B}_1, \*{B}_2$ be two Boolean algebras and $\rho: \*{B}_1 \xhookrightarrow{} \*{B}_2$ an embedding. We define the \emph{full CNf-twist-algebra over $\*{B}_1$ and $\*{B}_2$} as the algebra $\*{B}_1 \bowtie \*{B}_2$ with universe $B_1 \times B_2$ and operations
\begin{itemize}
    \setlength \itemsep{0pt}   
    \item[i.] $(x_1, x_2) \land (y_1, y_2) := (x_1 \land_{\*B_1} y_1, x_2 \lor_{\*B_2} y_2)$,
    \item[ii.] $\neg (x_1, x_2) := (\neg_{\*B_1} x_1, \neg_{\*B_2} x_2)$, 
    \item[iii.] $(x_1, x_2) \to (y_1, y_2) := (x_1 \to_{\*B_1} y_1, \rho(x_1) \to_{\*B_2} y_2)$, 
    \item[iv.] $\bot := (0_{\*B_1}, 0_{\*B_2}).$
    \end{itemize}
    Additionally, one can define $x \lor y := \neg (\neg x \land \neg y)$. A \emph{CNf-twist-algebra over $\*{B}_1$ and $\*{B}_2$} is any subalgebra $\*A\leqslant \*{B}_1 \bowtie \*{B}_2$ verifying that $\pi_1[A] = B_1$. 
\end{definition}
Note that in this case we have that $\top$ is definable by $\neg x \to (x \to x)$ and hence $\bot$ is too by letting $\bot := \neg \top$ (alternatively, we may set $\top := (0, 1) \to (1, 0)$ and $\bot := \neg \top$)\footnote{Again, we can make similar observations as those in Remark \ref{rmk:90deg}.}. Again, as expected, $\mathbf{CNf_4}$ is isomorphic to the CNf-twist-algebra over two copies of the two-element Boolean algebra (where $\rho$ is the identity).

\begin{remark}
\label{remark:booleanreducts}
Note that, in Definitions \ref{Fftwist}, \ref{DFftwist} and \ref{CNftwist} the structures defined have a Boolean reduct given by restricting to the operations $\land$ and $\neg$. 
\end{remark}

Consider the class of algebras: 

\begin{definition}
    A \emph{CNf-algebra} $(D; \land, \lor,  \to, \neg,0, 1, \top)$ is an expansion of a De Morgan algebra $(D; \land, \lor, \neg, 0, 1)$ with an operation $\to$ verifying: 
    \begin{enumerate}[(CNf1)]
        \item $\top \approx 1 \Rightarrow x \approx y$ (see Remark \ref{remark:quasieq}),
        \item $\top \land \bot = 0$,
        \item $\Diamond a = \Diamond b$ implies $\Box \neg(a \to \bot) = \Box \neg(b \to \bot)$,
        \item $\Diamond a = \Diamond b$ and $\Box a = \Box b$ imply that $a = b$,
        \item $\Box (a \to b) = \Box ((a \to \bot) \land b)$,
        \item $\Diamond (a \land \neg a) = \Box (a \lor \neg a) = \bot$, 
        \item $\Diamond (a \lor \neg a) = \Box (a \land \neg a) = 1$,
        \item $\Box \neg(a \to \bot) = \Box \neg(b \to \bot)$ implies that $\Diamond a = \Diamond b$,
        \item $\Diamond (a \to b) = \Diamond a \to \Diamond b$, 
        \item $\Box \neg ((a \land b) \to \bot) = \Box \neg ((a \to \bot) \land (b \to \bot))$,
        \item $\Box (\neg a \to \bot) =\Box \neg (a \to \bot)$,
        \item $\Box (a \to \bot \land \neg (b \to \bot)) = \Box \neg ((a \to b) \to \bot)$,
    \end{enumerate}
    where $\bot := \neg \top$, $\Diamond x := x \lor \bot$ and $\Box x := x \land \bot$. 
\end{definition}

Similarly as before, given a CNf-algebra $\*D$, define $\Diamond x := x \lor \bot$ and $\Box x := x \land \bot$. First, consider $\Diamond(\* D) := (\Diamond[D]; \land, \lor, \neg_\Diamond, \to, \bot, 1)$, where $\neg_\Diamond \Diamond a := \Diamond \neg a$. Then, let us define $\Box(\*D) := (\Box(D); \land_\Box, \lor_\Box, \neg_\Box, \to_\Box, 0_\Box, 1_\Box),$ where $\Box a \land_\Box \Box b := \Box (a \lor b)$, $\Box a \lor_\Box \Box b := \Box (a \land b)$, $\neg_\Box \Box a := \Box \neg a$ and $\Box a \to_\Box \Box b := \neg_\Box \Box a \lor_\Box \Box b = \Box (\neg a \land b)$, $0_\Box := \bot$ and $1_\Box := 0$. 

\begin{lemma}
\label{CNflemma}
   Given a CNf-algebra $\*D$, it holds that both $\Diamond(\*D)$ and $\Box(\*D)$ are Boolean algebras.
\end{lemma}
\begin{proof}
  We already know that $\Diamond(\* D)$ and $\Box(\*D)$ are De Morgan algebras. It is enough to see that they are complemented (Definition \ref{complemented}). But one can readily check that the following equalities hold by the definition above: 
  $\Diamond a \land \neg_\Diamond \Diamond a = \Diamond a \land \Diamond \neg a = \Diamond (a \land \neg a) = \bot$ and $\Diamond a \lor \neg_\Diamond \Diamond a = \Diamond a \lor \Diamond \neg a = \Diamond (a \lor \neg a) = 1$.
  Analogously, one can check that $\Box a \land_\Box \neg_\Box \Box a = \Box a \land_\Box \Box \neg a = \Box (a \lor \neg a) = \bot$ and $\Box a \lor_\Box \neg_\Box \Box a = \Box a \lor_\Box \Box \neg a = \Box (a \land \neg a) = 1$, as desired.
\end{proof}

Now, as before, we need to give an embedding between the two Boolean algebras under consideration: 

\begin{lemma}
\label{CNfembed}
    The map $\rho: \Diamond(\*D) \to \Box(\*D): \Diamond a \mapsto \Box \neg(a \to \bot)$ is an embedding of Boolean algebras.
\end{lemma}
\begin{proof}
We have to check the following conditions:
\begin{itemize}
    \item[--] $\rho$ is well defined. If $\Diamond a = \Diamond b$ then $\rho(\Diamond a) = \Box \neg(a \to \bot) = \Box \neg(b \to \bot) = \rho(\Diamond b)$ by definition.
    \item[--] $\rho$ is injective. If $\rho(\Diamond a) = \Box \neg(a \to \bot) = \Box \neg(b \to \bot) = \rho(\Diamond b)$ then $\Diamond a = \Diamond b$.
    \item[--] $\rho$ preserves $\land$. Indeed, $\rho(\Diamond (a \land b)) =\Box \neg ((a \land b) \to \bot) = \Box \neg ((a \to \bot) \land (b \to \bot)) = \Box (\neg(a \to \bot) \lor \neg (b \to \bot)) = \Box \neg(a \to \bot) \land_\Box \Box \neg(b \to \bot) = \rho(\Diamond a) \land_\Box \rho(\Diamond b)$.
    \item[--] $\rho$ preserves $\neg$. Indeed, $\rho(\neg \Diamond a) = \rho(\Diamond \neg a) = \Box \neg (\neg a \to \bot) = \neg_\Box \Box (\neg a \to \bot) = \neg_\Box \Box \neg (a \to \bot) = \neg_\Box \rho(\Diamond a)$. 
    \item[--] $\rho$ preserves $\to$. Indeed, $\rho(\Diamond a) \to_\Box \rho(\Diamond b) = \neg_\Box \rho(\Diamond a) \lor_\Box \rho(\Diamond b) = \neg_\Box(\Box \neg (a \to \bot)) \lor_\Box \Box \neg (b \to \bot) = \Box (a \to \bot \land \neg (b \to \bot))$ and this, by definition, is $\Box \neg ((a \to b) \to \bot) = \rho(\Diamond (a \to b)) = \rho(\Diamond a \to \Diamond b)$.
    \item[--] $\rho$ preserves extrema. Indeed, $\rho(\bot) = \rho(\Diamond \bot) = \Box \neg (\bot \to \bot) = \Box \neg \top = \Box \bot = \bot$. On the other hand, $\rho(1) = \rho(\Diamond 1) = \Box \neg (1 \to \bot) = \Box \neg \bot = 0$, as desired.
\end{itemize}  
\end{proof}

Now we are in conditions of proving: 

\begin{theorem}[CNf-twist representation]
\label{CNftwistrep}
    Let $\*D$ be a CNf-algebra. Then, the map $\iota: \*D \to  \Diamond (\*D)\times \Box(\*D): a \mapsto (\Diamond a, \Box a)$ witnesses an isomorphism between $\*D$ and the corresponding CNf-twist-structure.
\end{theorem}
\begin{proof}
We need to check the following requirements:
    \begin{itemize}
     \item[--] $\iota$ is injective. This means that $\Diamond a = \Diamond b$ and $\Box a = \Box b$ imply $a = b$. But this holds by definition.
    \item[--] $\iota$ preserves $\land$. We have $\iota(a \land b) = (\Diamond(a \land b), \Box (a \land b)) = (\Diamond a \land \Diamond b, \Box a \lor_\Box \Box b) = \iota(a) \land \iota(b),$ as we wanted.
\item[--] $\iota$ preserves $\neg$. Indeed, $\iota(\neg a) = (\Diamond \neg a, \Box \neg a) = (\neg_\Diamond \Diamond a, \neg_\Box \Box a) = \neg \iota(a)$.
    \item[--] $\iota$ preserves $\bot$. This holds because $\iota(\bot) = (\Diamond \bot, \Box \bot) = (\bot, \bot)$, as desired.
    \item[--] $\iota$ preserves $\to$. Here we have, on the one hand: 
    $$\iota(a \to b) = (\Diamond(a \to b), \Box (a \to b)) = (\Diamond a \to \Diamond b, \Box(a \to b)),$$
    and, on the other: 
    $$\iota(a) \to \iota(b) = (\Diamond a, \Box a) \to (\Diamond b, \Box b) = (\Diamond a \to \Diamond a, \rho(\Diamond a)\to_\Box \Box b),$$
    so we have to note that
    $$\rho(\Diamond a)\to_\Box \Box b = \Box\neg(a \to \bot) \to_\Box \Box b = \neg_\Box \Box \neg (a \to \bot) \lor_\Box \Box b = \Box ((a \to \bot) \land b)$$
    and that this equals $\Box(a \to b)$ by definition.
    \end{itemize}
\end{proof}

\paragraph{Cooper.} We define OLf as the logic induced by $(\mathbf{OL^f_4}, \{\top, \tv\})$, where $\mathbf{OL^f_4} := (A_4; \land_\mathsf{OL}^4, \lor_\mathsf{OL}^4, \to_\mathsf{OL}^4, \neg_\mathsf{f})$ and 
the table of $\neg_\mathsf{f}$ is the following:

\begin{displaymath}
\begin{array}{|c|c|}
 & \neg_\mathsf{f} \\ 
\hline 
\bv & \tv\\
\bot & \top\\
\top & \bot \\
\tv & \bv\\
\end{array}
\quad 
\begin{array}{|c|c c c c|}
\to_\mathsf{OL}^4 & \bv & \bot & \top & \tv\\ 
\hline 
\bv & \top & \top & \top & \top\\
\bot & \top & \top & \top & \top\\
\top & \bv & \bot & \top & \tv\\
\tv & \bv & \bot & \top & \tv \\
\end{array}
\quad 
\begin{array}{|c|c c c c|}
\land_\mathsf{OL}^4 & \bv & \bot & \top & \tv\\ 
\hline 
\bv & \bv & \bv& \bv & \bv\\
\bot & \bv & \bv & \bot & \bot\\
\top & \bv & \bot & \top & \tv\\
\tv & \bv & \bot& \tv & \tv \\
\end{array}
\quad 
\begin{array}{|c|c c c c|}
\lor_\mathsf{OL}^4 & \bv & \bot & \top & \tv\\ 
\hline 
\bv &  \bv &  \bot & \bv & \tv\\
\bot & \bot & \tv & \bot & \tv \\
\top & \bv & \bot &\top & \tv\\
\tv & \tv & \tv & \tv & \tv \\
\end{array}
\end{displaymath}

It is then clear that $\mathbf{OLf_4}$ is a particular example of the following twist construction by letting both factors to be the same and $\rho$ to be the identity:

\begin{definition}
\label{OLftwist}
     Let $_1$ and $B_2$ be two Boolean algebras and $\rho: B_1 \to B_2$ an embedding. An \emph{OLf-twist-algebra over $B_1$ and $B_2$} is an algebra $B_1 \bowtie B_2$ consisting of universe $B_1\times B_2$ and operations
\begin{itemize}
    \setlength \itemsep{0pt}   
    \item[i.] $(x_1, x_2) * (y_1, y_2) := (x_1 \land_{\*B_1} y_1, (\rho(x_1) \to_{\*B_2} y_2)\land_{\*B_2} (\rho(y_1) \to_{\*B_2} x_2)),$ 
    \item[ii.] $(x_1, x_2) \to (y_1, y_2) := (x_1 \to_{\*B_1} y_1, \rho(x_1) \to_{\*B_2} y_2),$ 
    \item[iii.] $\neg (x_1, x_2) := (\neg_{\*B} x_1, \neg_{\*B} x_2)$.
\end{itemize}
\end{definition}
The main desideratum here would be, again, a twist representation result as in \cite{Riv2x}. 
As we did in the OLg case, it seems natural to employ the comments from Remark \ref{OLalgCNalg}, although in this context its application is not so straightforward. Recall that we use $*$ for the operation $\land$ from the previous definition. First, one refines our previous definition in the following way: 

\begin{definition}
    \label{OLftwist2}
    Let $\*{B}_1, \*{B}_2$ be two Boolean algebras and $\rho: \*{B}_1 \xhookrightarrow{} \*{B}_2$ an embedding. We define the \emph{full OLf-twist-algebra over $\*{B}_1$ and $\*{B}_2$} as the full CNf-twist-algebra $\*{B}_1 \bowtie \*{B}_2$ over $\*B_1$ and $\*B_2$ expanded with the operation:
$$(x_1, x_2) * (y_1, y_2) := (x_1 \land_{\*B_1} y_1, (\rho(x_1) \to_{\*B_2} y_2)\land_{\*B_2} (\rho(y_1) \to_{\*B_2} x_2)).$$
    An \emph{OLf-twist-algebra over $\*{B}_1$ and $\*{B}_2$} is any subalgebra $\*A\leqslant \*{B}_1 \bowtie \*{B}_2$ verifying that $\pi_1[A] = B_1$. 
\end{definition}

Then, one needs to check that every OLf-twist-algebra in this sense belongs to the following class:

\begin{definition}
    \label{OLfalg}
    An \emph{OLf-algebra} $(D; \land, \lor,  \to, *, \neg, 0, 1, \top)$ is a CNf-algebra $(D; \land, \lor, \to,\neg,  0, 1, \top)$ expanded with an operation $*$ verifying that: 
    \begin{enumerate}[(OLf1)]
        \item $\Diamond (a * b) = \Diamond (a \land b)$,
        \item $\Box(a * b) = \Box (((a \to \bot) \land b) \lor ((b \to \bot) \land  a))$,
    \end{enumerate}
    where $\Diamond x := x \lor \bot$, $\Box x := x \land \bot$.
\end{definition}

Then, given a OLf-algebra $\*D$, we already know that we can obtain two Boolean algebras $\Diamond(\*D)$ and $\Box(\*D)$ as in Lemma \ref{CNflemma} and that, moreover, there is an embedding $\rho: \Diamond(\* D) \xhookrightarrow{} \Box(\*D)$. It is enough to check, then, that the map from Theorem \ref{CNftwistrep} preserves $*$:

\begin{theorem}[OLf-twist representation]
\label{OLftwistrep}
    Let $\*D$ be a OLf-algebra. Then, the map $\iota: \*D \to  \Diamond (\*D)\times \Box(\*D): a \mapsto (\Diamond a, \Box a)$ witnesses an isomorphism between $\*D$ and the corresponding OLf-twist-structure.
\end{theorem}
\begin{proof}
We only need to check that the definition above ensures the desired equalities. On one hand, we have that $\iota(a*b) = (\Diamond(a * b), \Box (a * b))$ while, on the other, we have 
$$\iota(a) * \iota(b) = (\Diamond a \land \Diamond b, (\rho(\Diamond a) \to_\Box \Box b)\land_\Box (\rho(\Diamond b) \to_\Box \Box a))).$$
But
$$(\rho(\Diamond a) \to_\Box \Box b)\land_\Box (\rho(\Diamond b) \to_\Box \Box a)) = (\Box \neg (a \to \bot)\to_\Box \Box b )\land_\Box (\Box \neg(b \to \bot) \to_\Box \Box a)$$
and this equals
\begin{align*}
\Box (\neg \neg (a \to \bot) \land b ) \land_\Box \Box(\neg \neg(b \to \bot) \land  a) 
& = \Box ((a \to \bot) \land b) \land_\Box  \Box((b \to \bot) \land  a) \\
& = \Box (((a \to \bot) \land b) \lor ((b \to \bot) \land  a)).
\end{align*}
So the result holds by definition.
\end{proof}

\section{Concluding remarks and future work}

The present paper has been concerned with presenting twist representation results for the three-valued logics of indicative conditionals considered, studying notable fragments of these, and with extending the three-valued setting to a four-valued following the intuitions provided by twist-algebras. 
Some central questions  have been only touched upon
in the preceding pages, and so stand in need of further development. We list a few below.

\paragraph{1. Analyzing the logical systems.} It would be desirable to study 
consequences of the twist representation results presented here. Along the lines of~\cite{Riv2x,wollic}, 
it would be interesting to establish algebraizability results and, more in general, obtaining
a deeper understanding of the algebraic semantics of these systems. 
Axiomatizability issues for our logics also will need to be addressed (as in the papers~\cite{wollic,tableaux}). 
Also following  the approach of~\cite{wollic}, one 
could study possible variations of the logics of interest.
Comparing the expressive power of the various logics (especially in the four-valued case)
also appears to be a promising line to be pursued.
Lastly, from a mathematical as well as a philosophical point of view, it would be interesting  to identify a
core set of formulas that could enable a classification of the logics 
according to whether they satisfy or not the corresponding principles
(the numerous examples from~\cite{Cooper1968} could, in this respect, provide a starting point).

\paragraph{2. Further philosophical analysis of the four-valued logics.} 
The fundamental philosophical issues of
three-valued logics of indicative conditionals have to do with 
(i) whether to designate the third value $\top$ and (ii) fixing the truth tables. The first, (i), is usually settled by appealing to either the \emph{non-falsity} or the \emph{information-preserving} account (as opposed to the classical truth-preserving notion), while (ii) depends on the informal reading of the connectives. 
Ultimately, one complements these answers by additional pragmatical considerations, such as designing a concrete set of desirable properties that the logic in question should meet (v.g. some axioms and rules, connexive principles of some kind, etc.).
%
In the present paper we have provided a few intuitions motivating that our new four-valued tables 
extend the three-valued ones in a natural way.
A deeper philosophical analysis motivating these logics needs still to be addressed, and would be worth exploring in future research.

\paragraph{3. Material collapse.} Many authors  have defended the material conditional as the most adequate one (see~\cite{egrecollapse} for a survey). The idea that any satisfactory theory of indicative conditionals will collapse, in some relevant sense, to the material conditional account can be traced back, at least, to Gibbard~\cite{gibbard}, who presents a general argument against the conception of $\to$ 
as a logical connective. Indeed, if $\to$ verifies 
minimal 
requirements 
usually postulated for any implication-like connective (notably the so-called Import-Export rule), then $\varphi \to \psi$ turns out to be inter-derivable with the material conditional $\neg \varphi \lor \psi$. This argument has been recently refined by Fitelson~\cite{fitelson} and further analyzed in the context of our three-valued logics in~\cite{egrecollapse}. 
Such a discussion could be extended to our new logics
in two alternative ways. One could argue that what Gibbard's argument achieves is only a weakening of a more general 
property that 
should be postulated 
by a strong proponent of material collapse, namely, self-extensionality. This property fails in the case of our logics~\cite{wollic}, and similar arguments could be made for the new four-valued logics. On the other hand, one could look at the Boolean behavior of some of these (Remark~\ref{remark:booleanreducts}) as a strong indication that, in some cases, we obtain systems not very far apart from classical logic.

\paragraph{4. Fragments of $\mathcal{NINE}$.} A possible objection
to 
the \emph{vacuously true} approach sketched above is the following: 
why are $\top$ and $\bot$ not comparable? For vacuous truth is, after all, a truth,
so
it should stand \emph{above} vacuous falsity in the logical order. 
One could thus explore alternative ways of ordering the truth-values. 
Suppose that, in our preceding sensor interpretation motivating $\bot$, 
we considered three instead of two possible states, namely: being \emph{weakly}, \emph{strongly} or \emph{not} detected. This framework would relate  our logics with the 
logic $\mathcal{NINE}$ of~\cite{damasio} and with the bilattice approach (essentially akin to twist-structures: see~\cite{ferguson,jungriv,Riv14}). 
For instance, we could identify the values $\bv, \bot, \top, \tv$ with the pairs $(none, strong)$, $(none, weak)$, $(weak, none)$, $(strong, none)$; note that $\land$ and $\neg$ can be defined as in our twist-algebras for the three-valued case. An advantage of this approach is that it may help 
account 
for the behavior of $\top$ in the original truth tables, for $\top$ behaves as a truly gappy/infectious value according to the truth table of $\to_\mathsf{DF}$, as a weak truth in the case of $\to_\mathsf{OL}$ and as a something in between 
according to the table of 
$\to_\mathsf{F}$. Identifying the truth-values with different fragments of $\mathcal{NINE}$ in each case might
provide further insight on this phenomenon. 



\paragraph{5. Other alternatives.} A further approach for extending the three-valued setting could consist in making use of non-deterministic matrices (see, e.g., \cite{nondetum,addingimp}), which would allow us to introduce $\bot$ as a free truth-value, letting the negation of $\bot$ be any of the four truth-values (this situation is philosophically different from the one mentioned in footnote \ref{footnote1}, for here the classical truth-values would behave as expected). Alternatively, one could  study linear orderings on $A_4$  together with the associated 
conjunction and disjunction, as we have briefly discussed. One could consider, for instance, 
the order $\bot < \bv < \tv < \top$ extending the one induced by $\land_\mathsf{OL}$: an interesting problem would then be that of introducing a meaningful implication operation that would interact nicely with this order.







\bibliographystyle{splncs04}
\bibliography{mybib}

@article{Riv2x,
	author = "Rivieccio, U.",
        journal = "Archive for Mathematical Logic",
	title = {The algebra of ordinary discourse. \uppercase{O}n the semantics of \uppercase{C}ooper’s logic},
	year = {2025}, 
        DOI={10.1007/s00153-024-00961-2}
}

@article{Cooper1968,
	author = {William S. Cooper},
	doi = {10.1080/00201746808601531},
	journal = {Inquiry: An Interdisciplinary Journal of Philosophy},
	number = {1-4},
	pages = {295--320},
	publisher = {Taylor \& Francis},
	title = {The \uppercase{P}ropositional \uppercase{L}ogic of \uppercase{O}rdinary \uppercase{D}iscourse},
	volume = {11},
	year = {1968}
}

@book{F16,
author = {Font, Josep Maria},
year = {2016},
month = {04},
title = {Abstract Algebraic Logic: An introductory textbook},
isbn = {978-1-84890-207-7},
publisher = {College Publications}
}

@book{BuSa00,
    author = "S. Burris and H. P. Sankappanavar",
    title = "A Course in Universal Algebra",
    publisher = "Springer",
    year = "2011"
}

@Book{ClDa98,
  Title                    = {Natural dualities for the working algebraist},
  Author                   = {D. M. Clark and B. A. Davey},
  Publisher                = {Cambridge University Press},
  Year                     = {1998},

  Address                  = {Cambridge},
  Series                   = {Cambridge Studies in Advanced Mathematics},
  Volume                   = {57},
  Filename                 = {ClDa98},
  ISBN                     = {0-521-45415-8},
  Mrnumber                 = {MR1663208 (2000d:18001)},
  Pages                    = {xii+356},
}

@article{egre20211,
author="{\'E}gr{\'e}, Paul
and Rossi, Lorenzo
and Sprenger, Jan",
title="De {F}inettian logics of indicative conditionals part {I}: trivalent semantics and validity",
journal="Journal of Philosophical Logic",
year="2021",
month="Apr",
day="01",
volume="50",
number="2",
pages="187--213",
issn="1573-0433",
doi="10.1007/s10992-020-09549-6",
}

@article{egre20212,
author="{\'E}gré, P.
and Rossi, L.
and Sprenger, J.",
title="{D}e {F}inettian logics of indicative conditionals part {II}: proof theory and algebraic semantics",
journal="Journal of Philosophical Logic",
year="2021",
month="Apr",
day="01",
volume="50",
number="2",
pages="215--247",
issn="1573-0433",
doi="10.1007/s10992-020-09572-7",
}

@article{FaOd23,
author="Fazio, Davide
and Odintsov, Sergei P.",
title="An Algebraic Investigation of the {C}onnexive {L}ogic {C}",
journal="Studia Logica",
year="2023",
month="Jun",
day="21",
issn="1572-8730",
doi="10.1007/s11225-023-10057-2"
}

@incollection{WanConnexive05,
	author = {Heinrich Wansing},
	booktitle = {Advances in Modal Logic},
	editor = {Marcus Kracht and Maarten de Rijke and Heinrich Wansing and Michael Zakharyaschev},
	pages = {367--383},
	publisher = {CSLI Publications},
	title = {Connexive Modal Logic},
	year = {2005}
}

@article{Um17,
author="Albuquerque, Hugo and
Přenosil, Adam and 
Rivieccio, Umberto",
title="{An Algebraic View of Super-Belnap Logics}",
journal="Studia Logica",
year="2017",
month="Dec",
volume="105",
number="6",
pages="1051--1086",
doi="10.1007/s11225-017-9739-7"
}

@article{Cig86,
    author = {Cignoli, R.},
    title = {The class of {K}leene algebras satisfying an interpolation property and {N}elson algebras},
    journal = "Algebra Universalis",
    volume="23",
    year = "1986",
    pages="262-292",
doi = "https://doi.org/10.1007/BF01230621"
}

@book{FJa09,
    author = "Font, J. M. and Jansana, R.",
    title = "{A General Algebraic Semantics for Sentential Logics}",
    publisher = "Springer",
    year = "2009" 
}

@article{Greati23,
    author = "Greati, V. and Marcelino, S. and Rivieccio, U.",
    title = "{Axiomatizing the Logic of Ordinary Discourse}",
    journal = "Proceeding of IPMU24 (to appear)",
    year = "May 5 2024"
}

@article{Fa86,
    author = "Farrell, R. J.",
    title = "{Implication and Presupposition}",
    journal = "Notre Dame Journal of Formal Logic, Volume 27, Number 1",
    year = "January 1986" 
}

@article{Um20,
    author = "Rivieccio, U.",
    title = "Representation of {D}e {M}organ and ({S}emi-){K}leene Lattices",
    journal = "Soft Computing, 24, pp. 8685–8716 ",
    year = "2020",
    doi = "https://doi.org/10.1007/s00500-020-04885-w"
}

@article{Ca08,
    author = "Cantwell, J.",
    title = "{The Logic of Conditional Negation}",
    journal = "Notre Dame Journal of Formal Logic 49(3): 245-260",
    year = "2008", 
    doi = "10.1215/00294527-2008-010"
}

@article{OmWa20,
    author = "Omori, H. and Wansing, H.",
    title = " An extension of connexive logic C",
    journal = "In: Olivetti, N., Verbrugge, R., Negri, S., Sandu, G. (eds.) Advances in Modal Logic, 13, pp. 503–522. College Publications",
    year = "2005"
}

@article{deFin,
    author = "de Finetti, B.",
    title = "La logique de la probabilité",
    journal = "Actes du Congrés International de Philosophie Scientifique",
    year = "1936"
}

@article{priest,
    author = "Priest, G.",
    title = "{The Logic of Paradox}",
    journal = "Journal of Philosophical Logic, Vol. 8, No. 1, pp. 219-241",
    year = "Jan., 1979"
}

@article{Fa79,
    author = "Robert J. Farrell",
    title = "Material implication, confirmation, and counterfactuals",
    journal = "Notre Dame Journal of Formal Logic 20 (2), pp. 383-394",
    year = "1979"
}

@article{tableaux,
    author = "Vitor Greati and Sergio Marcelino and Miguel {Muñoz Pérez} and Umberto Rivieccio",
    title = "Analytic calculi for logics of indicative conditionals",
    journal = "In: Pozzato, G. L., Uustalu, T. (eds) Automated Reasoning with Analytic Tableaux and Related Methods. TABLEAUX 2025. Lecture Notes in Computer Science(), vol 15980. Springer, Cham.",
    year = "2026"
}

@article{wollic,
    author = "Miguel {Muñoz Pérez} and Umberto Rivieccio",
    title = "{Indicative Conditionals: Some Algebraic Considerations}",
    journal = "In Kozen, D. \& de Queiroz, R. (eds) Logic, Language, Information, and Computation, WoLLIC 2025, Lecture Notes in Computer Science, vol 15942, Springer",
    year = "2025"
}

@book{humberstone,
    author = "Lloyd Humberstone",
    title = "The connectives",
    publisher = "MIT Press" ,
    year = "2011"
}

@article{avron,
    author = "Ofer Arieli and Arnon Avron",
    title = "{Four-Valued Paradefinite Logics}",
    journal = "Studia Logica 105, pp. 1087–1122, Springer",
    year = "2017"
}

@article{nondetum,
    author = "C Caleiro and S Marcelino and U Rivieccio",
    title = "{Some More Theorems on Structural Entailment Relations and Non-deterministic Semantics}",
    journal = "In Malinowski, J., Palczewski, R. (eds) Janusz Czelakowski on Logical Consequence, Outstanding Contributions to Logic, vol 27. Springer",
    year = "2024"
}

@article{addingimp,
    author = "V Greati and S Marcelino and J Marcos and U Rivieccio",
    title = "Adding an implication to logics of perfect paradefinite
algebras",
    journal = "Mathematical Structures in Computer Science 34, pp. 1138–1183",
    year = "2024"
}

@article{jungriv,
    author = "A Jung and U Rivieccio",
    title = "{Priestley Duality for Bilattices}",
    journal = "Studia Logica 100 (1-2):223-252",
    year = "2012"
}

@article{khoo,
    author = "J Khoo",
    title = "On indicative and subjunctive conditionals",
    journal = "Philosophers' Imprint, Vol. 15, No. 32",
    year = "2015"
}

@article{JaRi21,
    author = "R Jansana and U Rivieccio",
    title = "{Quasi-Nelson algebras and fragments}",
    journal = "Mathematical Structures in Computer Science, 31, 2021, pp. 257–285",
    year = "2021"
}

@article{szmuc,
    author = "B {Da Ré} and F Pailos and D Szmuc",
    title = "Theories of truth based on four-valued infectious logics",
    journal = "Logic Journal of the IGPL, Volume 28, Issue 5, Pages 712–746",
    year = "2020"
}

@article{ferguson,
    author = "T M Ferguson",
    title = "{Cut-Down Operations on Bilattices}",
    journal = "45th International Symposium on Multiple-Valued Logic",
    year = "2015"
}

@article{damasio,
    author = "C V Damásio and L M Pereira",
    title = "A model theory for paraconsistent logic programming",
    journal = "In: Pinto-Ferreira, C., Mamede, N.J. (eds) Progress in Artificial Intelligence. EPIA 1995. Lecture Notes in Computer Science, vol 990. Springer, Berlin, Heidelberg",
    year = "1995"
}

@article{sobocinski,
    author = "B Sobocinski",
    title = "Axiomatization of a partial system of three-valued calculus of propositions",
    journal = "The Journal of Computing Systems, vol. 1, pp. 23-55",
    year = "1952"
}

@article{Pa72,
    author = "R Z Parks",
    title = "A note on {R}-mingle and {S}obocinski’s three-valued logic",
    journal = " Notre Dame Journal of Formal Logic 13.2, pp. 227-228",
    year = "1972"
}

@article{egrecollapse,
    author = "P {\'E}gré and L Rossi and J Sprenger",
    title = "Gibbardian {C}ollapse and {T}rivalent {C}onditionals",
    journal = "In: Kaufmann, S., Over, D.E., Sharma, G. (eds) Conditionals. Palgrave Studies in Pragmatics, Language and Cognition. Palgrave Macmillan, Cham.",
    year = "2023"
}

@article{gibbard,
    author = "A Gibbard",
    title =  "Two {R}ecent {T}heories of {C}onditionals",
    journal = "In: Harper, W.L., Stalnaker, R., Pearce, G. (eds) IFS. The University of Western Ontario Series in Philosophy of Science, vol 15. Springer, Dordrecht",
    year = "1980"
}

@article{fitelson,
    author = "B Fitelson",
    title = "Gibbard’s Collapse Theorem for the Indicative Conditional: An Axiomatic Approach",
    journal = "In: Bonacina, M.P., Stickel, M.E. (eds) Automated Reasoning and Mathematics. Lecture Notes in Computer Science(), vol 7788. Springer, Berlin, Heidelberg",
    year = "2013"
}

@article{BlRa08,
    author = "W J Blok and G Raftery",
    title = "{Assertionally Equivalent Quasivarieties}",
    journal = "International Journal of Algebra and Computation, Vol. 18, No. 04, pp. 589-681",
    year = "2008"
}

@article{BlRa04,
    author = "W J Blok and J G Raftery",
    title = "{Fragments of R-Mingle}",
    journal = "Studia Logica 78, pp.59–106",
    year = "2004"
}

@article{Co21,
    author = "M E Coniglio and T G G{ó}mez & M Figallo",
    title = "{Some model-theoretic results on the 3-valued paraconsistent first-order logic QCiore}",
    journal = "The Review of Symbolic Logic, Vol. 14, No. 1",
    year = "2021"
}

@article{CaMa00,
    author = "W A Carnielli and J Marcos and S {de Amo}",
    title = "{Formal inconsistency and evolutionary databases}",
    journal = "Logic and Logical Philosophy 8 (2):115-152",
    year = "2000"
}

@article{Riv14,
    author = "U Rivieccio",
    title = "{Implicative twist-structures}",
    journal = "{Algebra Universalis, 71, pp. 155–186}",
    year = "2014" 
}

@article{Ra06,
    author = "J G Raftery",
    title = "{The Equational Definability of Truth Predicates}",
    journal = "Reports on Mathematical Logic 41, pp. 95-149",
    year = "2006"
}

@article{celani08,
    author = "S Celani and L M Cabrer",
    title = "{Topological duality for Tarski algebras}",
    journal = "{Algebra universalis 58 (1), 73-94}",
    year = "2008"
}

@article{egreprob,
    author = "{{\'E}gré, P. and Rossi, L. and Sprenger, J.}",
    title = "Probability for Trivalent Conditionals",
    journal = "Preprint",
    year = "2024"
}

@article{sprenger,
    author = "J Sprenger",
    title = "{The Conditional in Three-Valued Logic}",
    journal = "{In: Égré, P. \& Rossi, L. (eds) Handbook of Trivalent Logics. Cambridge, Massachusetts: The MIT Press (forthcoming)}",
    year = "2024"
}
\end{document}